\documentclass{article}
\usepackage{amsmath,amsthm,amsfonts,amscd,eucal}

\numberwithin{equation}{section}

\def\ca{{\mathcal A}}
\def\cb{{\mathcal B}}
\def\cc{{\mathcal C}}

\def\cf{{\mathcal F}}

\def\ch{{\mathcal H}}
\def\cai{{\mathcal I}}

\def\ck{{\mathcal K}}

\def\cam{{\mathcal M}}

\def\cu{{\mathcal U}}

\def\bc{{\mathbb C}}

\def\bn{{\mathbb N}}
\def\br{{\mathbb R}}

\def\a{\alpha}
\def\b{\beta}

\def\eps{\varepsilon}

\def\l{\lambda}       
\def\m{\mu}
\def\n{\nu}
\def\x{\xi}
\def\r{\rho}
\def\s{\sigma}
    
\def\t{\tau}

\def\f{\varphi}

\def\o{\omega}        \def\Q{\Omega}\def\O{\Omega}

\def\imply{\Rightarrow}

\def\ov{\overline}
\def\itm#1{\item{$(#1)$}}

\newcommand{\id}{\text{\rm id}}
\newcommand{\dist}{\text{\rm dist}}

\DeclareMathOperator{\vol}{vol}

\def\wt{\widetilde}
\def\cst{\text{C$^{*}$-algebra}}

\def\ov{\overline}

\newcommand{\diam}{\text{diam}}
\def\lsp{Lip-space}
\newcommand{\diag}{\text{diag}}

\def\nl#1{\|#1\|_{L}}

\def\tn#1{|\hskip-1pt|\hskip-1pt|#1|\hskip-1pt|\hskip-1pt|}

\newtheorem{Thm}{Theorem}[section]
\newtheorem{Cor}[Thm]{Corollary}
\newtheorem{Prop}[Thm]{Proposition}
\newtheorem{Lemma}[Thm]{Lemma}

\theoremstyle{definition}
\newtheorem{Dfn}[Thm]{Definition}

\theoremstyle{remark}
\newtheorem{rem}[Thm]{Remark} 
\newtheorem{ack}{Acknowledgement} 
\begin{document}

\title{The problem of completeness for Gromov-Hausdorff metrics on
C$^{*}$-algebras
\thanks{Supported in part by GNAMPA, MIUR, EU Research Training
Network in Quantum Spaces - Noncommutative Geometry.
E-mail: {\tt guido@mat.uniroma2.it, isola@mat.uniroma2.it}}} 

\author{Daniele Guido, Tommaso Isola
\\
Dipartimento di Matematica, \\ Universit\`a di Roma ``Tor
Vergata'', \\  I--00133 Roma, Italy.
}

\date{}
\maketitle
\markboth{Completeness for Gromov-Hausdorff metrics 
on C$^{*}$-algebras}
{Completeness for Gromov-Hausdorff metrics 
on C$^{*}$-algebras}

\begin{abstract}
    It is proved that the family of equivalence classes of Lip-normed
    C$^{*}$-algebras introduced by M. Rieffel, up to complete order
    isomorphisms preserving the Lip-seminorm, is not complete w.r.t.
    the matricial quantum Gromov-Hausdorff distance introduced by D.
    Kerr.  This is shown by exhibiting a Cauchy sequence whose limit,
    which always exists as an operator system, is not completely order
    isomorphic to any \cst.
    
    Conditions ensuring the existence of a C$^{*}$-structure on the
    limit are considered, making use of the notion of ultraproduct. 
    More precisely, a necessary and sufficient condition is given for
    the existence, on the limiting operator system, of a
    C$^{*}$-product structure inherited from the approximating
    \cst{s}.  Such condition can be considered as a generalisation of
    the $f$-Leibniz conditions introduced by Kerr and Li. 
    Furthermore, it is shown that our condition is not necessary for
    the existence of a C$^{*}$-structure {\it tout court}, namely
    there are cases in which the limit is a \cst, but the
    C$^{*}$-structure is not inherited.
\end{abstract}

\section{Introduction}

In this paper we study the problem of completeness of some spaces
w.r.t the matricial quantum Gromov-Hausdorff metrics introduced by
Kerr \cite{Kerr}, showing that the space of equivalence classes of
\cst{s} with Lipschitz seminorms is not complete.  

As is known, Rieffel introduced and studied, in a series of papers
\cite{Rieffel98,Rieffel99,Rieffel00,Rie4,Rie5,Rie6}, the notion of
compact quantum metric space, and generalised the Gromov-Hausdorff
distance to the quantum case.  The main tool is a seminorm $L$ on the
``quantum'' functions, which plays the role of the Lipschitz seminorm
for functions on a compact metric space.  The requirements can be
formalised as follows: $L$ should vanish exactly on the multiples of
the identity element, and should induce on the states (positive
normalised functionals) the weak$^{*}$ topology.  It is not
restrictive to assume that $L$ is also lower-semicontinuous in norm,
and we shall always assume it in this paper.  A space endowed with
such a seminorm is called Lip-normed.

Roughly speaking, the quantum Gromov-Hausdorff distance between two
C$^{*}$-algebras corresponds to the Gromov-Hausdorff distance between
the corresponding state spaces, endowed with the Monge-Kantorovitch
metric induced by the Lipschitz seminorms (for the precise definition
see eq. (\ref{qGHdist})).  However one easily realises that while two
abelian Lip-normed C$^{*}$-algebras having zero quantum
Gromov-Hausdorff distance are isomorphic, the same is not true for
noncommutative C$^{*}$-algebras.  The structure which is preserved by
the quantum Gromov-Hausdorff distance is indeed that of order-unit
space.

In fact Rieffel proved that the quantum Gromov-Hausdorff distance is
indeed a distance on equivalence classes of order-unit spaces, showed
that the space of equivalence classes is complete, and gave also
conditions for compactness.

As mentioned above, when \cst{s} are concerned, there are non
isomorphic Lip-normed \cst{s} which have zero quantum Gromov-Hausdorff
distance.  In order to cope with this problem, Kerr \cite{Kerr}
introduced matricial quantum Gromov-Hausdorff distances $\dist_{p}$. 
When $p$ is finite, $\dist_{p}$ measures the distance between $p\times
p$-valued state spaces, and $\dist_{\infty}$ corresponds to the
supremum over all $p$.

He showed that, when $p\geq 2$, $\dist_{p}$ vanishes if and only if the
Lip-normed \cst{s} are $^{*}$-isomorphic.  The question of completeness for
the space of equivalence classes remained open.  However Kerr and Li introduced
a family of conditions \cite{Kerr,Li2}, depending on a function $f$,
related with the Leibniz property for the Lipschitz norms, showing
that if all the Lip-normed \cst{s} of a Cauchy sequence satisfy the
same $f$-Leibniz condition they converge to a \cst{} (which satisfies
the same $f$-Leibniz condition).

The main purpose of this paper is then to solve the completeness
problem, indeed we exhibit a Cauchy sequence of Lip-normed \cst{s}
which does not converge to a \cst.

Following Kerr, a natural setting for matricial quantum
Gromov-Hausdorff distance is that of Lip-normed operator systems.  He
showed that the distance $\dist_{\infty}$ between two Lip-normed
operator systems vanishes if and only if they are completely order
isomorphic and the Lipschitz norm is preserved by the isomorphism.

It is possible to show that the space of equivalence classes of
operator systems with Lipschitz seminorms, endowed with
$\dist_{\infty}$, is complete (we do so in Theorem \ref{thm:2.12}, by
making use of ultraproducts, and it was also proved independently by
Kerr and Li \cite{KL} with different techniques).  Therefore the
problem of completeness for Lip-normed \cst{s} w.r.t. the 
$\dist_{\infty}$
metric can be rephrased as the problem of the closure of the
Lip-normed \cst{s} inside the family of equivalence classes of
Lip-normed operator systems.

Given a $\dist_{\infty}$-Cauchy sequence of Lip-normed \cst{s}, the
limit always exists as an operator system $S$, and the question
becomes to determine whether such $S$ admits a C$^{*}$-structure (is
completely order isomorphic to a \cst).  This property can be
reformulated with the aid of the notion of injective envelope
$\cai(S)$ of an operator system $S$ due to Hamana \cite{Ha}.  The
operator system $S$ embeds canonically in $\cai(S)$, and the latter
admits a unique C$^{*}$-product.  The existence of a C$^{*}$-structure
on $S$ is equivalent to the fact that $S$ is a subalgebra of
$\cai(S)$.  We use this technique in subsection \ref{ex1} to show that
the limit of a suitable sequence of \cst{s} with Lipschitz seminorms
is not a \cst.

As mentioned above, an important tool in our analysis is the notion of
ultraproduct of Banach spaces, and in particular a tailored version
for Banach spaces with lower semicontinuous Lip\-schitz norm, which we
call Lip-ultraproduct.

We show that under a condition of uniform compactness on sequences of
spaces, the Lip-ultraproduct is a Banach space with lower
semicontinuous Lip\-schitz norm, and inherits some of the structures
from the approximating spaces, in particular those of order-unit space
and of operator system.  Furthermore we show that Cauchy sequences are
uniformly compact and the Lip-ultraproduct is indeed the limit.

Let us mention here that the representative of a quantum
Gromov-Hausdorff limit constructed via ultraproducts is directly
endowed with a lower semicontinuous Lip-seminorm.

The C$^{*}$-structure is not inherited in general by the
Lip-ultraproduct, however for any given free ultrafilter the
Lip-ultraproduct is always a closed linear subspace of the
ultraproduct, and the latter is a \cst.  Therefore there are cases in
which the limit inherits a C$^{*}$-structure, namely when the
Lip-ultraproduct is a subalgebra of the ultraproduct (for a suitable
free ultrafilter $\cu$).  This is a sufficient condition for the limit to
be a \cst, but is not necessary, namely there are cases in which the
limit is a \cst{} but the C$^{*}$-structure is not inherited, cf. 
subsection \ref{ex2}.

Moreover we can completely characterize the Cauchy sequences for which
the limit inherits a C$^{*}$-structure in terms of a function
$\eps(r)$, $r\in[0,+\infty)$, associated with any Lip-normed \cst,
measuring how far is the set of Lipschitz elements from being an
algebra.  If we have a Cauchy sequence with functions $\eps_{n}(r)$,
the limit inherits the C$^{*}$-structure if and only if, for a
suitable subsequence $n_{k}$, $\limsup_{k} \eps_{n_{k}}(r) \to 0$ for
$r\to\infty$ (cf.  Corollary \ref{Cor:inherited}).  Therefore such
condition is a maximal generalisation of the $f$-Leibniz condition of
Kerr and Li.  The fact that it is indeed more general is illustrated
in subsection \ref{ex3}.

Let us also mention that, with the aid of the function $\eps(r)$ and
of the results on inherited C$^{*}$-structure, we can easily
manufacture a new distance on the family of equivalence classes of
Lip-normed \cst{s}, for which completeness holds, cf.  Corollary
\ref{Cor:completeness}.  The convergence condition under this new
distance is clearly stronger than the convergence condition w.r.t.
$\dist_{\infty}$, as shown by Example 2 in subsection \ref{ex2}.  However
this stronger convergence condition seems to be more natural when
\cst{s} are concerned, because in this case the C$^{*}$-structure is
always inherited, namely the product on the limit is the limit of the
approximating products, cf.  equation (\ref{inherprod}).

As mentioned above, subsection \ref{ex1} is devoted to the
construction of examples of non converging Cauchy sequences w.r.t. the
matricial quantum Gromov Hausdorff distance.  When the limit (as an
operator system) does not inherit a C$^{*}$-structure, as a
Lip-ultraproduct it is described by a subspace, which is not closed
w.r.t. the product, of a \cst{} (the ultraproduct).  For the examples
considered in subsection \ref{ex1} we show that the product structure
given by the immersion in the ultraproduct is the same as the product
structure given by the immersion in the injective envelope, thus
showing that the limit is not a \cst.

Indeed the examples considered in subsection \ref{ex1} depend on a
\cst{} $\cb$, and we show that for any $\cb$ we get a sequence
$\ca_{n}$ which is Cauchy w.r.t. $\dist_{p}$, $p\in\bn\cup\{\infty\}$.  In
the particular case in which $\cb=\bc I$, the sequence $\ca_{n}$
consists of the constant algebra $M_{2}(\bc)$ of $2\times 2$ matrices, and
it is easy to show that the limit is not even positively isomorphic to
a \cst{} (cf. Remark \ref{rem:UHF}).  This shows that the family of equivalence
classes of Lip-normed \cst{s} is not complete w.r.t. $\dist_{p}$,
$p\geq2$.  However, if we confine our attention to the case
$\ca_{n}=M_{2}(\bc)$, one may argue that we have simply chosen the wrong
distance.

Let us recall that when Rieffel introduced the quantum
Gromov-Hausdorff distance, he had to generalise to the quantum setting
a distance involving spaces of points, or extremal states.  Since for
\cst{s} extremal states may be not closed, and even dense, as in the
UHF case, he decided to consider a distance involving all states. 
However, when $\ca_{n}=M_{2}(\bc)$, the replacement of the quantum
Gromov-Hausdorff distance with a distance involving only extremal
states, like the distance $\dist^{e}_{q}$ considered by Rieffel in
\cite{Rieffel00} after Proposition 4.9, would destroy the
counterexample, since the sequence is no longer Cauchy w.r.t. such
distance.

This is the reason why we consider also non-trivial $\cb$: when the
C$^{*}$-algebra $\cb$ is UHF, we get a sequence made of a constant UHF
algebra (with different Lip-norms), for which pure states are dense,
hence matricial quantum Gromov-Hausdorff distances are the only
reasonable choices.  Of course in this case the proof that the limit
is not completely order isomorphic to a \cst{} is more difficult,
requiring the notion of injective envelope of Hamana \cite{Ha}.

We conclude by mentioning a result for ultraproducts which may have an
interest of its own.  The dual of an ultraproduct is larger in general
than the ultraproduct of the duals, the equality being attained only
under a strong uniform convexity property of the sequence, which is
never satisfied for infinite-dimensional \cst{s}.  For the
Lip-ultraproduct however, if the sequence is uniformly compact, any
element in the dual can be realised as an element in the ultraproduct
of the dual spaces, namely the compactness condition of the Lipschitz
seminorms allows one to construct a more manageable ultraproduct,
whose dual is made of equivalence classes of sequences of functionals.

This suggests the interpretation of the Lip-ultraproduct as the
quantum (dualised) analogue of the ultralimit of compact metric
spaces.  As in the classical case, an ultralimit is a limit only if a
uniform compactness condition is satisfied.

\section{Order-unit spaces}

This section is mainly devoted to the introduction of the 
Lip-ultraproduct and the study of its properties.

In order to clarify some features of the construction, we introduce 
the notion of Lip-space.

Let us recall (see \cite{Rieffel98}, Thm.  1.9) that a lower
semicontinuous Lipschitz seminorm $L$ on a complete order-unit space
can be characterised, besides the vanishing exactly on the multiples
of the identity, by the fact that the elements whose norm and
Lipschitz seminorm are bounded by a constant, form a compact set in
norm.  Indeed by introducing the norm
$\|x\|_{L}=\max\{L(x),\frac1{R}\|x\|\}$, where $R$ may be taken as
half of the diameter of the state space w.r.t. the Lipschitz distance,
the compactness property may be reformulated as the fact that the
$\|\cdot\|_{L}$-balls are norm compact, and the Lipschitz seminorm can
be recovered as $L(x)=\inf_{\lambda}\|x-\l I\|_{L}$.  Therefore, in
contrast with the standard terminology, we shall reserve the term
Lip-norm for $\|\cdot\|_{L}$, and shall call $L$ a Lip-seminorm.

The observations above suggest the definition of a Lip-space as a 
Banach space with  an extra norm $\|\cdot\|_{L}$ (finite on a dense 
subspace) such that the $\|\cdot\|_{L}$-balls are compact.

\subsection{{\lsp}s}

\begin{Dfn}\label{def:2.1}
    We call {\lsp}  a triple $(X,\|\cdot\|,\nl{\cdot})$
    where 
    \item{$(i)$} $(X,\|\cdot\|)$ is a Banach space, 
    \item{$(ii)$} $\nl{\cdot} : X\to [0,
    +\infty]$ is finite on a dense vector subspace $X_{0}$ where it is a 
    norm,
    \item{$(iii)$}  the unit ball w.r.t. $\nl{\cdot}$, $\{x\in 
    X:\nl{x}\leq 1\}$, is compact in $(X,\|\cdot\|)$.
\end{Dfn}

We call {\it radius} of the {\lsp} $(X,\|\cdot\|,\nl{\cdot})$, and
denote it by $R$, the maximum of $\|\cdot\|$ on the unit ball w.r.t.
$\|\cdot\|_{L}$, hence
\begin{equation}
    \|x\|\leq R\nl{x},\qquad x\in X.
\end{equation}
As we shall see it is the analogue of the radius of a Lip-normed order
unit space introduced by Rieffel at the end of Section 2 in
\cite{Rieffel00}.

\begin{Prop}\label{Prop1.1a}
    Let $(X,\|\cdot\|,\nl{\cdot})$ be a Lip-space, $e\in 
    X_{0}\setminus\{0\}$, and set $L(x) := \inf_{\l\in\br} \| x-\l
    e\|_{L}$. Then $L$ is a lower semicontinuous densely defined 
    seminorm and $L(x)=0 \iff x=\l e$ for some $\l\in\br$.
\end{Prop}
\begin{proof}
    Indeed, it is easy to
    prove that $L$ is a seminorm, and that $L(\l e)=0,\l\in\br$. 
    Moreover, as $\|\cdot\|_{L}$ is lower semicontinuous, because of 
    Definition \ref{def:2.1} $(iii)$, and $\|x-\l e\|_{L}\geq |\l| 
    \|e\|_{L} - \|x\|_{L}\to\infty, |\l|\to\infty$, we obtain $L(x)= 
    \min_{\l\in\br} \|x-\l e\|_{L}$.   
    Therefore, if $L(x)=0$, then there is
    $\l_{0}\in\br$ s.t. $\| x-\l_{0}e\|_{L}=0$, so that $x=\l_{0}e$. 
    
    Finally, if $x, x_{n}\in X$, $\| x_{n}-x\|\to0$,
    then, $L(x) \leq \liminf_{n\to\infty} L(x_{n})$. Indeed, 
    passing possibly to a subsequence, we may assume $\{L(x_{n})\}$ 
    converges. Let, for all $n\in\bn$, $\l_{n}\in\br$ be
    s.t. $\| x_{n}-\l_{n} e\|_{L} = L(x_{n})$. Then $\{\| x_{n}-\l_{n} 
    e\|_{L}\}$ is bounded; 
    so by Definition \ref{def:2.1} $(iii)$, there are 
    $\{n_{k}\}\subset\bn$, 
    $a\in X$ s.t. $\| x_{n_{k}}-\l_{n_{k}} e - a\|\to0$. 
    Therefore there is $\l_{0}\in\br$ s.t. $\l_{n_{k}}\to\l_{0}$, and 
    $a=x-\l_{0}e$. Hence
    \begin{align*}
	L(x) & \leq \|x-\l_{0}e\|_{L} \leq \liminf_{k\to\infty}
	\|x_{n_{k}}-\l_{n_{k}} e\|_{L} \\
	& = \lim_{k\to\infty} L(x_{n_{k}}) = \lim_{n\to\infty}
	L(x_{n}),
    \end{align*}
    where the second inequality follows from Definition 
    \ref{def:2.1} $(iii)$.
\end{proof}

\begin{Prop}\label{prop:1.2}
    Let $(X,\|\cdot\|,\nl{\cdot})$ be a {\lsp}. Then 
    the dual norm
    $$
    \nl{\x}':=\max_{x\in X}\frac{|\langle\x,x\rangle|}{\nl{x}}
    $$
    induces the weak$^{*}$ topology on the bounded subsets of $X'$, 
    the Banach space dual of $(X,\|\cdot\|)$.
        
    The constant $R$ is equal to the radius, in the $\nl{\cdot}'$ norm,
    of the unit ball of $(X',\|\cdot\|)$.
\end{Prop}

\begin{proof}
    First observe that $\nl{\cdot}'$, which is obviously a seminorm, 
    is indeed a norm. In fact, if $\nl{\x}'=0$, then $\x$ vanishes 
    on $X_{0}$, which is dense, i.e. $\x=0$.
    
    Now we consider the identity map $\iota$ from the closed unit ball
    $B'_{1}$ of $X'$ endowed with the weak$^{*}$ topology to the same
    set endowed with the distance induced by $\nl{\cdot}'$.  Given
    $r>0$, we consider a $r/2$-net $\{x_{i}:i=1,\dots,n\}$ in $\{x\in
    X:\nl{x}\leq 1\}$.  Then, if $\|\x\|'\leq1$,
    $$
    |\langle\x,x\rangle| \leq
    \max_{i=1,\dots,n}|\langle\x,x_{i}\rangle| + r/2.
    $$
    Therefore, the weak$^{*}$ open set in $B'_{1}$
    $$
    U=\{\|\x\|'\leq1:\max_{i=1,\dots,n}|\langle\x,x_{i}\rangle|<r/2\},
    $$
    is contained in the $\nl{\cdot}'$ open set in $B'_{1}$
    $$
    V=\{\|\x\|'\leq1:\nl{\x}'<r\},
    $$
    showing that $\iota$ is continuous. Since the domain is compact 
    and the range is Hausdorff, $\iota$ is indeed a homeomorphism.
    
    Finally, the radius of the unit ball of $X'$ in the $\nl{\cdot}'$
    norm is
    $$
    \sup_{\|\x\|'\leq1}\nl{\x}'=\sup_{\x\ne0,x\ne0}
    \frac{|\langle\x,x\rangle|}{\nl{x}\|\x\|'}
    =\sup_{x\ne0}\frac{\|x\|}{\nl{x}}
    \sup_{\x\ne0}\frac{|\langle\x,x\rangle|}{\|\x\|'\|x\|}
    =R.
    $$
\end{proof}

\begin{Dfn}
    A family $\cf$ of {\lsp}s is called uniform if for all
    $\eps>0$ there is $n_{\eps}\in\bn$ such that, for any
    $(X,\|\cdot\|,\nl{\cdot})$ in $\cf$, $\{x\in X:\nl{x}\leq
    1\}$ can be covered by $n_{\eps}$ $\|\cdot\|$-balls of radius $\eps$.
\end{Dfn}

\begin{Lemma}\label{unifR}
    If $\cf$ is a uniform family of {\lsp}s, there is $R>0$ such that
    $\|x\|\leq R\nl{x}$ for any $(X,\|\cdot\|,\nl{\cdot})$ in $\cf$,
    $x\in X$.
\end{Lemma}
    
\begin{proof}
    Let $(X,\|\cdot\|,\nl{\cdot})$ be a {\lsp} such that $\{x\in
    X:\nl{x}\leq 1\}$ can be covered by $n$ balls of radius $1$, and
    let $x_{0}\in X$, $\nl{x_{0}}=1$.  Since the set $\{tx_{0}:
    t\in[0,1]\}$ is contained in $\{x\in X:\nl{x}\leq 1\}$, it is
    covered by at most $n$ balls of radius $1$, hence its length is
    majorised by $2n$, i.e. $R\leq2n$.
\end{proof}

\begin{Lemma}\label{eqnorms}
    Let $(V,\|\cdot\|)$ be an $n$-dimensional normed space.  Then the
    ball of radius $R$ can be covered by $(2R/\eps)^{n}$ balls of
    radius $\eps$.
\end{Lemma}

\begin{proof}
    Let us recall that, denoting by $n_{\eps}(\O)$ the minimum number
    of balls of radius $\eps$ covering $\O$, and by $\n_{\eps}(\O)$
    the maximum number of disjoint balls of radius $\eps$ contained in
    $\O$, one gets $n_{\eps}(\O)\leq\n_{\eps/2}(\O)$ (cf.  e.g.
    \cite{GuIs6}, Lemma 1.3).  Then, denoting by $\vol$ the Lebesgue
    measure and by $B_{r}$ the ball of radius $r$ w.r.t. the given
    norm, we get $\vol(B_{R}) \geq \n_{\eps}(B_{R}) \vol(B_{\eps})$,
    and $\vol(B_{R})= (R/\eps)^{n} \vol(B_{\eps})$, hence
    $n_{\eps}(B_{R})\leq (2R/\eps)^{n}$.
\end{proof}

\begin{Prop}\label{equivUniform}
    A family $\cf$ of {\lsp}s is uniform $\Leftrightarrow$ there
    exists a constant $R$ as in Lemma \ref{unifR}, and
    $\forall\eps>0$, there is $N_{\eps}\in\bn$ such that any {\lsp}
    $X$ in $\cf$ has a subspace $V$ of dimension not greater than
    $N_{\eps}$ such that $\{x\in V : \nl{x}\leq 1\}$ is $\eps$-dense
    in $\{x\in X:\nl{x}\leq 1\}$.
\end{Prop}

\begin{proof}
    $(\Rightarrow)$ The constant $R$ exists
    by Lemma \ref{unifR}; choose a covering of $\{x\in
    X:\nl{x}\leq 1\}$ by $n_{\eps}$ $\|\cdot\|$-balls of radius $\eps$ and
    consider the vector space $V$ generated by their centres.  Its
    dimension is clearly majorised by $n_{\eps}$.
    \\
    $(\Leftarrow)$ Take $\eps\leq1$.  The elements in $\{x\in V :
    \nl{x}\leq 1\}$ are contained in $\{x\in V : \|x\|\leq R\}$, hence
    any covering of the $R$-normic ball of $V$ with balls of radius
    $\eps$ gives a covering of the Lip-norm unit ball in $X$ with balls of
    radius $2\eps$.  By Lemma \ref{eqnorms}, one can realise the
    former covering with $(2R/\eps)^{N_{\eps}}$ balls, hence the
    implication is proved.
\end{proof}

\subsection{Ultraproducts}

Given a sequence $(X_{n},\|\cdot\|,\nl{\cdot})$ of {\lsp}s, we may
consider the Banach space $\ell^{\infty}(X_{n})$ of norm-bounded
sequences $x_{n}\in X_{n}$ with the sup-norm.  As is known \cite{Si},
if $\cu$ is a free ultrafilter on $\bn$, the ultraproduct
$\ell^{\infty}(X_{n},\cu)$ is defined as the quotient of
$\ell^{\infty}(X_{n})$ w.r.t. the subspace of sequences such that
$\lim_{\cu}\|x_{n}\|=0$.  We denote by $\pi_{\cu}$ the projection from
$\ell^{\infty}(X_{n})$ onto $\ell^{\infty}(X_{n},\cu)$.

\begin{Dfn}
    Given a sequence $(X_{n},\|\cdot\|,\nl{\cdot})$ of {\lsp}s, we
    call Lip-ultraproduct, and denote it by
    $\ell^{\infty}_{L}(X_{n},\cu)$, or simply by $X_{\cu}$, the image
    under $\pi_{\cu}$ of $\ell^{\infty}_{L}(X_{n})$, the norm closure
    of the space of bounded sequences for which
    $\nl{\{x_{n}\}}:=\sup_{\bn}\nl{x_{n}}<+\infty$.
\end{Dfn}

The quotient norm $\|\cdot\|_{\cu}$ of the equivalence class 
$x_{\cu}$ of a sequence $x_{n}$ is defined as
$$
\|x_{\cu}\|_{\cu}=
\inf_{[y_{n}]= x_{\cu}}
\sup_{n}\|y_{n}\|,
$$
hence $\|x_{\cu}\|_{\cu}=\lim_{\cu}\|x_{n}\|$ (\cite{AkKh}, Chap. 2 Prop. 
2.3).

Analogously, the quotient norm $\|\cdot\|_{L,\cu}$ of  
$x_{\cu}$ is defined as
\begin{equation}\label{quotientLip}
    \|x_{\cu}\|_{L,\cu}=
    \inf_{[y_{n}]= x_{\cu}}
    \sup_{n}\|y_{n}\|_{L}.
\end{equation}
This implies that $\|x_{\cu}\|_{L,\cu}\leq\lim_{\cu}\|x_{n}\|_{L}$,
in fact for any $\eps>0$ there exists an element $U$ of the
free ultrafilter such that, for any $n\in U$, $\|x_{n}\|_{L} \leq
\lim_{\cu}\|x_{m}\|_{L}+\eps$.  Then we may define $y_{n}=x_{n}$
for $n\in U$ and $y_{n}=0$ for $n\not\in U$.  Since
$[y_{n}]= [x_{n}]$, the result follows.

\begin{Lemma}\label{Lem:1.9}
    The infimum in $(\ref{quotientLip})$ is indeed a minimum.
\end{Lemma}
\begin{proof}
    Given $x_{\cu}\in X_{\cu}$, we may choose sequences $x^{k}_{n}$
    realising it and such that $\|x^{k}_{n}\|_{L}\leq
    \|x_{\cu}\|_{L,\cu}(1+\frac1k)$.  It is also not restrictive to
    ask that all the vectors $x^{k}_{n}$ have norm bounded by
    $2\|x_{\cu}\|$.  Then we set
    \begin{align*}
	V_{k}&=\{n\geq k:\|x_{n}^{j}-x_{n}^{i}\|\leq\frac1i,i\leq j\leq k\},\\
	V_{0}&=\bn,
    \end{align*}
    and observe that $V_{k}\in\cu$, $V_{k+1}\subseteq V_{k}$, and
    $\displaystyle{ \bigcup_{k\geq0}V_{k} \setminus V_{k+1}=\bn}$.
    Then we define 
    $$
    \tilde{x}_{n}=\frac{k}{k+1} x^{k}_{n},\quad n\in V_{k}\setminus
    V_{k+1},
    $$
    implying $\|\tilde{x}_{n}\|_{L}\leq \|x_{\cu}\|_{L,\cu}$.
    Now we show that $\tilde{x}_{\cu} = x_{\cu}$.  Indeed, if $n\in
    V_{i}$, $\exists k\geq i$ s.t. $n\in V_{k}\setminus V_{k+1}$, hence
    $$
    \|\tilde{x}_{n}-x^{i}_{n}\|\leq
    \|\tilde{x}_{n}-x^{k}_{n}\| + \|x^{k}_{n}-x^{i}_{n}\|\leq
    \frac1{k+1}\|x^{k}_{n}\|+\frac1i\leq(2\|x_{\cu}\|_{\cu}+1)\frac1i.
    $$
    Since $n$ is eventually in $V_{i}$ w.r.t. $\cu$, we get
    $$
    \|\tilde{x}_{\cu}-x_{\cu}\|=
    \lim_{\cu} \|\tilde{x}_{n}-x^{i}_{n}\|\leq
    \left(2\|x_{\cu}\|_{\cu}+1\right)\frac1i.
    $$
    By the arbitrariness of $i$ we get the result.
\end{proof}

Choosing $\tilde{x}_{n}$ as in the proof above, we get
$$\|x_{\cu}\|_{L,\cu} = \lim_{\cu}\|\tilde{x}_{n}\|_{L}
=\sup_{n}\|\tilde{x}_{n}\|_{L}.
$$
In particular we obtain that, for any element 
$x\in\ell^{\infty}_{L}(X_{n},\cu)$,
\begin{equation}\label{minLip}
    \|x\|_{L,\cu}=\min_{[x_{n}]=x}\lim_{\cu}\|x_{n}\|_{L}.
\end{equation}

\begin{Prop}\label{restrUltra}
    Given a uniform sequence $(X_{n},\|\cdot\|,\nl{\cdot})$ of
    {\lsp}s, the Lip-ultraproduct $\ell^{\infty}_{L}(X_{n},\cu)$,
    endowed with the quotient norms $\|\cdot\|_{\cu}$,
    $\|\cdot\|_{L,\cu}$, is a {\lsp}.  Moreover, the radius $R$ for
    $\ell^{\infty}_{L}(X_{n},\cu)$ is equal to $\lim_{\cu}R_{n}$,
    where $R_{n}$ is the radius of $X_{n}$.
\end{Prop}

\begin{proof}
    Let us show that the closed Lip-norm unit ball in $X_{\cu}$ is
    totally bounded in norm.  Indeed, since given $\eps>0$, the closed
    Lip-norm unit ball in $X_{n}$ is covered by $n_{\eps}$ balls of
    radius $\eps$, we may choose points $x_{n,1},\dots,
    x_{n,n_{\eps}}$ in $X_{n}$ such that the closed Lip-norm ball of
    radius 2 in $X_{n}$ is covered by
    $$
    \bigcup_{i=1}^{n_{\eps}}B(x_{n,i},2\eps).
    $$
    Now, given any sequence $\{x_{n}\}_{n\in\bn}$, $x_{n}\in X_{n}$,
    $\|x_{n}\|_{L}\leq2$, we get
    $$
    \min_{i=1,\dots,n_{\eps}}\|x_{\cu}-x_{\cu,i}\|_{\cu}
    =\min_{i=1,\dots,n_{\eps}}\lim_{\cu}\|x_{n}-x_{n,i}\|
    =\lim_{\cu}\min_{i=1,\dots,n_{\eps}}\|x_{n}-x_{n,i}\|\leq 2 \eps.
    $$
    Since any $x_{\cu}$ s.t $\|x_{\cu}\|_{L,\cu}\leq1$ can be realized
    with a sequence $x_{n}$ such that $\|x_{n}\|_{L}\leq 2$, we get
    that the closed Lip-norm unit ball in $X_{\cu}$ is covered by
    $n_{\eps}$ balls of radius $3\eps$.
    
    Then we show that the closed Lip-norm unit ball in $X_{\cu}$ is norm
    complete, hence compact.  In fact let $\{x^{k}\}_{k\in\bn}$ be a
    Cauchy sequence of elements of $X_{\cu}$,
    $\|x^{k}\|_{L,\cu}\leq1$, and, according to the argument above,
    realize them via sequences $x^{k}_{n}$ such that
    $\|x^{k}_{n}\|_{L}\leq 1$.  Let us choose a diagonal sequence as
    follows.
    
    Set $\eps_{k}=\sup_{i,j\geq k}\|x^{i}-x^{j}\|$, and observe 
    that $\eps_{k}\to 0$. Then we consider the sets $V_{k}\subset\bn$ 
    defined as
    \begin{align*}
	V_{k}&=\{n\geq k:\|x^{j}_{n}-x^{i}_{n}\|\leq 2\eps_{i}, i\leq j\leq 
	k\}\\
	V_{0}&=\bn,
    \end{align*}
    and observe that $V_{k+1}\subseteq V_{k}$, $\displaystyle{
    \bigcup_{k\geq0}V_{k} \setminus V_{k+1}=\bn}$, and since
    $\lim_{\cu}\|x^{j}_{n}-x^{i}_{n}\|=\|x^{j}-x^{i}\|\leq \eps_{i}$, then
    $V_{k}\in\cu$.  Now we define the diagonal sequence as
    $$
    \tilde{x}_{n}=x^{k}_{n},\qquad n\in V_{k}\setminus V_{k+1}.
    $$
    Then, when $n\in V_{i}$, and $k\geq i$ satisfies $n\in
    V_{k}\setminus V_{k+1}$, we have $\|\tilde{x}_{n} - x^{i}_{n}\| =
    \|x^{k}_{n} - x^{i}_{n}\| \leq 2\eps_{i}$. Since $n$ is eventually 
    in $V_{i}$ w.r.t. $\cu$, $\|\tilde{x}_{\cu} -
    x^{i}\| = \lim_{\cu} \|\tilde{x}_{n} - x^{i}_{n}\| \leq 2\eps_{i}$,
    namely $\tilde{x}_{\cu}$ is the limit of the sequence $x^{k}$. 
    Therefore $\|\tilde{x}_{\cu}\|_{L,\cu} \leq \lim_{\cu}
    \|\tilde{x}_{n}\|_{L} \leq1$, i.e. the result.
    
    Finally we compute the constant $R$.  Let $x_{n}\in X_{n}$ be s.t.
    $\|x_{n}\|_{L}=1$, $\|x_{n}\|=R_{n}$, and consider the element
    $x_{\cu}\in X_{\cu}$.  As observed above, $\|x_{\cu}\|_{L}\leq1$
    and $\|x_{\cu}\|=\lim_{\cu}R_{n}$, implying
    $R\geq\lim_{\cu}R_{n}$.  Now, given $y_{\cu}\in X_{\cu}$ with $\|
    y_{\cu} \|_{L} \leq 1$, realise it via a sequence $y_{n}$ s.t. $\|
    y_{n} \|_{L} \leq 1$.  By definition, $\|y_{n}\|\leq R_{n}$,
    therefore
    $$
    \|y_{\cu}\|=\lim_{\cu}\|y_{n}\|\leq\lim_{\cu}R_{n},
    $$
    implying $R\leq\lim_{\cu}R_{n}$. The thesis follows.
\end{proof}
 
 The rest of this subsection is devoted to the study of the relation 
 between $\ell^{\infty}_{L}(X_{n},\cu)'$ and $\ell^{\infty}_{L}(X'_{n},\cu)$.

 \begin{Prop}\label{Prop1.6}
     Let $\{\s_{n}\in X'_{n}\}$ be uniformly bounded, and denote by
     $\s_{\cu}(x_{\cu}) := \lim_{\cu} \s_{n}(x_{n})$,
     $[x_{n}]=x_{\cu}\in \ell^{\infty}_{L}(X_{n},\cu)$.  Then
     $\s_{\cu}$ is well-defined, $\s_{\cu}\in
     \ell^{\infty}_{L}(X_{n},\cu)'$, and
     $$
     \|\s_{\cu}\|_{L,\cu}' = \lim_{\cu} \|\s_{n}\|_{L}'.
     $$
 \end{Prop}
 \begin{proof}
     Let $M>0$ be s.t. $\|\s_{n}\|'\leq M$, $n\in\bn$. We first prove that 
     $\s_{\cu}$ is well defined and bounded. Indeed, if 
     $[x_{n}']=[x_{n}]\in \ell^{\infty}_{L}(X_{n},\cu)$, then
     $\lim_{\cu} |\s_{n}(x_{n}')-\s_{n}(x_{n})| \leq 
     M\lim_{\cu}\|x_{n}'-x_{n}\| =0$. Moreover 
     $|\s_{\cu}(x_{\cu})|\leq M\lim_{\cu}\|x_{n}\| = M\|x_{\cu}\|$, so 
     that $\|\s_{\cu}\|_{\cu}'\leq M$.
     Finally
     \begin{align*}
	 \lim_{\cu} \|\s_{n}\|_{L}' & = \lim_{\cu}
	 \sup_{x_{n} \in X_{n}}
	 \frac{|\s_{n}(x_{n})|}{ \|x_{n}\|_{L} } =
	 \sup_{\{x_{n}\}\in \ell^{\infty}_{L}(X_{n})} \lim_{\cu}
	 \frac{|\s_{n}(x_{n})|}{ \|x_{n}\|_{L} } \\
	 & = \sup_{\{x_{n}\}\in \ell^{\infty}_{L}(X_{n})}
	 \frac{\lim_{\cu}|\s_{n}(x_{n})|}{\lim_{\cu} \|x_{n}\|_{L} }
	 = \sup_{x_{\cu}\in \ell^{\infty}_{L}(X_{n},\cu)} 
	 \sup_{[x_{n}]=x_{\cu}}
	 \frac{|\s_{\cu}(x_{\cu})|}{ \lim_{\cu} \|x_{n}\|_{L} } \\
	 & = \sup_{x_{\cu}\in \ell^{\infty}_{L}(X_{n},\cu)}
	 \frac{|\s_{\cu}(x_{\cu})|}{ \|x_{\cu}\|_{L,\cu} } =
	 \|\s_{\cu}\|'_{L,\cu},
     \end{align*}
     where in the last but one equality we used (\ref{minLip}).  Note
     also that, in that equality, the set of allowed elements in the
     supremum on the right is tacitly assumed not to contain
     $x_{\cu}=0$, while the set of allowed elements in the supremum on
     the left might also contain $x_{\cu}=0$, since in some examples
     one may find sequences $\{x_{n}\}$ such that $[x_{n}]=0$ but
     $\lim_{\cu}\|x_{n}\|_{L}>0$.  However for such sequences the
     numerator $|\s_{\cu}(x_{\cu})|$ is zero, therefore the supremum
     does not change.
 \end{proof}

\begin{Thm}\label{thm:dualUltraprod}
    Given a uniform sequence $(X_{n},\|\cdot\|,\nl{\cdot})$ of
    {\lsp}s, the ultraproduct $\ell^{\infty}(X_{n}',\cu)$ of the dual
    spaces projects on the dual $\ell^{\infty}_{L}(X_{n},\cu)'$ of
    the Lip-ultraproduct.  Moreover, given a uniformly bounded
    sequence $\s_{n}$ of elements in $X_{n}'$, the element $\s_{\cu}$
    in $\ell^{\infty}(X_{n}',\cu)$ gives the null functional on
    $\ell^{\infty}_{L}(X_{n},\cu)$ if and only if
    $\lim_{\cu}\|\s_{n}\|_{L}'=0$.
\end{Thm}

\begin{proof}
    We already observed that an element in $\ell^{\infty}(X_{n}',\cu)$
    gives rise to a functional on $X_{\cu}$, and since
    $\|\cdot\|_{L,\cu}$ is a norm on $(X_{\cu})'$, the last statement
    follows from Proposition \ref{Prop1.6}.
    
    It is known (\cite{Si}, Lemma 1, p.  77), and easy to show, that
    the pairing between $\ell^{\infty}(X_{n}')$ and
    $\ell^{\infty}(X_{n})$ given by $\langle\{\f_{n}\},
    \{x_{n}\}\rangle=\lim_{\cu}\f_{n}(x_{n})$ gives rise to a pairing
    between $\ell^{\infty}(X_{n}',\cu)$ and
    $\ell^{\infty}(X_{n},\cu)$, hence to an isometric map
    $\ell^{\infty}(X_{n}',\cu)\to \ell^{\infty}(X_{n},\cu)'$.  We are
    interested in the contraction $\pi:\ell^{\infty}(X_{n}',\cu)\to
    \ell^{\infty}_{L}(X_{n},\cu)'$ obtained by composing the previous
    isometric map with the quotient map from
    $\ell^{\infty}(X_{n},\cu)'$ to $\ell^{\infty}_{L}(X_{n},\cu)'$. 
    We have to show that $\pi$ is surjective.
    
    Given $\eps>0$, let us choose the subspaces $V_{n}\subset X_{n}$
    as in Proposition \ref{equivUniform}; we may also assume that all
    vectors in $V_{n}$ have finite Lip-norm, hence the $V_{n}$ form a
    uniform sequence of Lip-spaces with dimension bounded by
    $N_{\eps}$.  Clearly the Lip-ultraproduct $V_{\cu}$ can be seen as
    a subspace of $X_{\cu}$ of dimension at most $N_{\eps}$, and the
    Lip-norm unit ball of $V_{\cu}$ is $\eps$-dense in the Lip-norm
    unit ball of $X_{\cu}$.
    
    Since the $V_{n}$ have uniformly bounded dimension, $\ell^{\infty}
    (V_{n},\cu)' \equiv \ell^{\infty}(V_{n}',\cu)$ (cf.  \cite{Si},
    Theorem 2, p. 78).  Now take any norm-one element $\f \in
    (X_{\cu})'$, restrict it to $V_{\cu}$ and then extend it by
    Hahn-Banach theorem to an element $\tilde{\f}$ acting on
    $\ell^{\infty}(V_{n},\cu)$.  $\tilde{\f}$ can then be identified
    with an element of $\ell^{\infty}(V_{n}',\cu)$, namely we may find
    elements $\tilde{\f}_{n} \in V_{n}'$ such that $\tilde{\f} =
    [\tilde{\f}_{n}]$, $\|\tilde{\f}\|' = \lim_{\cu}
    \|\tilde{\f}_{n}\|'$.  Extend then $\tilde{\f}_{n}$ to an element
    $\f'_{n}\in X'_{n}$, and set $\f' := [\f'_{n}]\in
    \ell^{\infty}(X'_{n},\cu)$.  Clearly $\|\f'\|\leq1$, hence 
    $\pi(\f')\equiv \f'_{\cu}$,
    has norm less than 1, and observe that, by construction, $\f$ and
    $\pi(\f')$ coincide on $V_{\cu}$.
  
    For any element $x$ in $X_{\cu}$ with $\|x\|_{L}\leq1$ we may find
    $x_{\eps}\in V_{\cu}$ such that $\|x-x_{\eps}\|\leq\eps$,
    therefore
    \begin{equation*}
	|\f(x)-\pi(\f')(x)| \leq |\f(x) - \f(x_{\eps})| +
	|\pi(\f')(x_{\eps}) - \pi(\f')(x)| \leq 2\eps,
    \end{equation*}
    As a consequence,
    \begin{equation*}
	\|\pi(\f')-\f\|'_{L,\cu}=\sup_{\|x\|_{L}\leq1}|\f(x)-\pi(\f')(x)|
	\leq 2\eps.
    \end{equation*}
    Choosing $\eps=1/2k$, we may then construct sequences
    $\f^{k}_{n}\in X_{n}'$ such that $\|\f^{k}_{n}\|\leq1$ and,
    setting $\f^{k}=[\f^{k}_{n}]$, $\|\pi(\f^{k})-\f\|'_{L,\cu}\leq 1/k$. 
    Then we construct a diagonal sequence as in the proof of
    Proposition \ref{restrUltra}.
    
    Consider the sets $V_{k}\subset\bn$ defined as
    $$
    V_{k}=\{n\geq k:\|\f^{j}_{n}-\f^{i}_{n}\|'_{L}\leq \frac3i, i\leq j\leq k\},
    $$
    and observe that the $V_{k}$'s are non-increasing and belong to
    $\cu$.  Now we define the diagonal sequence as
    $$
    \tilde{\f}_{n}=\f^{k}_{n},\qquad n\in V_{k}\setminus V_{k+1}.
    $$
    Then, when $n\in V_{k}$, and $k'\geq k$ satisfies $n\in
    V_{k'}\setminus V_{k'+1}$, we have $\|\tilde{\f}_{n} -
    \f^{k}_{n}\|'_{L} = \|\f^{k'}_{n} - \f^{k}_{n}\|'_{L} \leq 3/k$
    hence, denoting by $\tilde{\f}$ the element in
    $\ell^{\infty}(X',\cu)$ corresponding to the sequence
    $\tilde{\f}_{n}$, we get $\|\pi(\tilde{\f}) - \pi(\f^{k})\|'_{L,\cu} =
    \lim_{\cu} \|\tilde{\f}_{n} - \f^{k}_{n}\|'_{L} \leq 3/k$, hence
    $\|\pi(\tilde{\f}) - \f\|'_{L,\cu}\leq\|\pi(\tilde{\f}) -
    \pi(\f^{k})\|'_{L,\cu} +\|\f - \pi(\f^{k})\|'_{L,\cu}\leq 4/k$.  By the
    arbitrariness of $k$ we get $\pi(\tilde{\f})=\f$.
\end{proof}

\subsection{Order-unit spaces}\label{subsec:1.3}

In this subsection the results obtained thus far are used to prove
that a Cauchy sequence of Lip-normed order-unit spaces converges to
the Lip-ultraproduct for any free ultrafilter, thereby providing a
different proof of a result already established by Rieffel
\cite{Rieffel00}, namely the completeness of the space of equivalence
classes of Lip-normed order-unit spaces w.r.t. the quantum Gromov
Hausdorff distance.  In section \ref{sec:3} the same approach,
suitably modified, will prove the completeness of the space of
equivalence classes of Lip-normed operator systems w.r.t.
$d_{\infty}$, a result recently proved by Kerr and Li (though with
different methods).  

We recall now the definition of order-unit space, referring to
\cite{Alfsen} for more details.

An \emph{order-unit space} is a real partially ordered vector space,
$X$, with a distinguished element $e$ (the order unit) satisfying: \\
1) (Order unit property) For each $a\in X$ there is an $r\in \br$ such
that $a\le re$;
\\
2) (Archimedean property) For $a\in X$, $a\le re$ for all $r>0\imply a\le 0$.

On an order-unit space $(X, e)$, we can define a norm as
\begin{eqnarray*}
\| a\|=\inf \{r\in \br: -re\le a\le re\}.
\end{eqnarray*}
Then $X$ becomes a normed vector space and we can consider its
dual, $X'$, consisting of the bounded linear functionals, equipped
with the dual norm $\| \cdot \|'$.

By a {\it state} of an order-unit space $(X, e)$, we mean a $\o \in
X'$ such that $\o(e)=\| \o\|' =1$.  States are automatically
positive.  Denote the set of all states of $X$ by $S(X)$.  It is a
compact convex subset of $X'$ under the weak$^*-$topology.  Kadison's
basic representation theorem \cite{Alfsen} says that the natural
pairing between $X$ and $S(X)$ induces an isometric order isomorphism
of $X$ onto a dense subspace of the space $Af_{\br}(S(X))$ of all
affine $\br-$valued continuous functions on $S(X)$, equipped with the
supremum norm and the usual order on functions. We denote by 
$\hat{a}(\o) :=\o(a),\ \o\in S(X)$, the affine function corresponding 
to $a\in X$.

For an order-unit space $(X, e)$, we say that a densely defined 
seminorm $L$ is a \emph{Lip-seminorm} (cf.  \cite[Definition
2.1]{Rieffel00}, where it is called Lip-norm) if:\\
1) For $a\in X$, we have $L(a)=0$ if and only if $a\in \br e$.\\
2) The topology on $S(X)$ induced by the metric $\rho_L$
\begin{equation} 
    \rho_{L}(\o_{1}, \o_{2})=\sup_{L(a)\le 1} |\o_{1}(a)-\o_{2}(a)|
\end{equation}
is the weak$^*-$topology.

We shall call Lip-normed order-unit space a complete order-unit space
endowed with a lower semicontinuous Lip-seminorm.

Let us recall that the {\it radius} $R$ of a Lip normed order-unit
space is defined as half of the diameter of $(S(X),\r_{L})$.  We now
endow $X$ with the norm
$$
\|a\|_{L}:= \max\{ \frac{\|a\|}{R}, L(a)\}.
$$
In the following Proposition we prove, for the sake of completeness,
some results which are needed in the sequel, even though some of them
are already known.
    
\begin{Prop}\label{prop:1.13} Let $(X,e,L)$ be a Lip-normed order-unit
space.  Then \itm{i} $\|a\|_{0} := \inf_{\l\in\br} \|a-\l e\| =
\frac12(\max\hat{a}-\min\hat{a})$,
    \itm{ii} $\|a\|_{L,0}:= \inf_{\l\in\br} \|a-\l e\|_{L} = L(a)= 
    \min_{\l\in\br}\|a-\l e\|_{L}$,
    \itm{iii} $R=\sup_{L(a)\neq 0} \frac{\|a\|_{0}}{L(a)}$, 
    \itm{iv} $\displaystyle{
    R
    = \sup_{\o\in S(\ca)} \|\o\|_{L}'
    = \sup_{\f\in\ca^{*},\|\f\|=1} \|\f\|_{L}' 
    = \sup_{a\neq 0} \frac{\|a\|}{\|a\|_{L}}}$.
\end{Prop}
\begin{proof}
$(i)$     
     \begin{align*}
	\|a\|_{0} & = \inf_{\l\in\br} \|a-\l e\| 
	 = \inf_{\l\in\br} \sup_{\o\in S(X)} |\hat{a}(\o)-\l| \\
	& = \inf_{\l\in\br} \max\{ |\max\hat{a}-\l|, 
	|\min\hat{a}-\l| \} 
	 = \frac{\max\hat{a}-\min\hat{a}}{2}.
    \end{align*}
$(ii)$ 
\begin{align*}
    \|a\|_{L,0}:
    &= \inf_{\l\in\br} \|a-\l e\|_{L} =
    \inf_{\l\in\br} \max\left\{\frac{\|a-\l e\|}{R}, L(a-\l e)\right\} \\
    &=\max\left\{ \inf_{\l\in\br} \frac{\|a-\l e\|}{R}, L(a)\right\}
    = \max\left\{ \frac{\|a\|_{0}}{R}, L(a)\right\} = L(a). 
\end{align*}
Because $\|\cdot\|_{L}$ is lower semicontinuous, the last
equality follows.\\
$(iii)$    
\begin{align*}
    \diam{\, S(X)} &:= \sup_{\o_{1},\o_{2}\in S(X)} \rho_{L}(\o_{1},\o_{2}) 
    = \sup_{\o_{1},\o_{2}\in S(X)} \sup_{L(a)\neq0} 
    \frac{|\o_{1}(a)-\o_{2}(a)|}{L(a)} \\
    & = \sup_{L(a)\neq0} 
    \frac{\displaystyle{\sup_{\o_{1},\o_{2}\in S(X)} |\o_{1}(a)-\o_{2}(a)|}}
    {L(a)} 
    = \sup_{L(a)\neq0} \frac{\displaystyle{\max_{\o\in S(X)} \hat{a}(\o) -
    \min_{\o\in S(X)} \hat{a}(\o)}}{L(a)} \\
    & = \sup_{L(a)\neq0} \frac{2\|a\|_{0}}{L(a)}.\\
\end{align*}
$(iv)$ Let us observe that $\|e\|_{L}=R^{-1}$, therefore

$$
R=\sup_{\o\in S(X)}\frac{\o(e)}{\|e\|_{L}} \leq
\sup_{\o\in S(X)}\|\o\|'_{L}.$$
Conversely,
$$
\sup_{\o\in S(X)}\|\o\|'_{L} = \sup_{\o\in S(X)}\sup_{a\in 
X}\frac{\o(a)}{\|a\|_{L}} \leq \sup_{\o\in S(X)}\sup_{a\in 
X}R\frac{\o(a)}{\|a\|}=R,
$$
proving the first equality.

As for the second, let $\f\in X'$, $\|\f\|=1$.  Then, from
\cite{Alfsen} II.1.14, there are $\r,\s\in S(X)$, $\l,\m\in[0,1]$,
$\l+\m=1$, s.t. $\f=\l\r-\m\s$.  Therefore
\begin{align*}
    \nl{\f}' &= \sup_{\nl{a}\leq1} |\f(a)| 
    \leq \sup_{\nl{a}\leq1} \left(\l|\r(a)| +\m |\s(a)|\right) \\
    & \leq \l \nl{\r}' +\m \nl{\s}' \\
    & = \sup_{\o\in S(X)} \nl{\o}',
\end{align*}
giving the result.

Finally,
$$
\sup_{\o\in S(X)}\|\o\|'_{L} 
= \sup_{\o\in S(X)}\sup_{a\in X}\frac{\o(a)}{\|a\|_{L}}
= \sup_{a\in X}\sup_{\o\in S(X)}\frac{\o(a)}{\|a\|_{L}}
=  \sup_{a\in X}\frac{\|a\|}{\|a\|_{L}}.
$$
\end{proof}

\begin{Thm}\label{thm:1.12}
    Let $(X,e,L)$ be a Lip-normed order-unit space of radius $R$, and
    define $\|a\|_{L} := \max\{ \frac{\|a\|}{R}, L(a)\}$, $a\in X$.
    Then $(X,\|\cdot\|,\|\cdot\|_{L})$ becomes a \lsp{} whose radius 
    as a Lip-space coincides with its radius as a Lip-normed order 
    unit space.
\end{Thm}
\begin{proof}
    As $\{a\in X: \nl{a}\leq 1\} = \{a\in X: \|a\|\leq R, L(a)\leq
    1\}$ is compact (\cite{Rieffel98}, Thm.  1.9), we get a 
    Lip-space. The equality between the radii follows from Proposition 
    \ref{prop:1.13} $(iv)$.
\end{proof}

 \begin{Prop}\label{prop:1.15}
     Let $\{(X_{n},e_{n})\}$ be complete order-unit spaces, $\cu$ a
     free ultrafilter.  Then the ultraproduct
     $(\ell^{\infty}(X_{n},\cu),e_{\cu})$ is a complete order-unit space.
 \end{Prop}
 \begin{proof}
     Let us recall that $\ell^{\infty}(X_{n},\cu) :=
     \ell^{\infty}(X_{n})/\cai_{\cu}$, where $\ell^{\infty}(X_{n})
     := \{ \{a_{n}\} : a_{n}\in X_{n}, \| \{a_{n}\} \| :=
     \sup_{n}\|a_{n}\|<\infty\}$, and $\cai_{\cu} := \{ \{a_{n}\}\in
     \ell^{\infty}(X_{n}) : \lim_{\cu} \| a_{n} \| = 0\}$.
         
     Observe that $\cai_{\cu}$ is a positively generated order ideal,
     because for any $\{a_{n}\}\in\cai_{\cu}$, there are
     $a_{n+},a_{n-}\in X_{n,+}$ s.t. $a_{n}=a_{n+}-a_{n-}$ and
     $\|a_{n\pm}\|\leq\|a_{n}\|$, see \cite{Alfsen} II.1.2. 
     
     Therefore, by \cite{Alfsen} II.1.6, we only have to check the
     Archimedean property for $\ell^{\infty}(X_{n},\cu)$.  Assume
     $a_{\cu}\leq \eps e_{\cu}$, for all $\eps>0$.  Then $\eps e_{\cu}
     -a_{\cu}\geq 0$, for all $\eps>0$, that is there is
     $U_{\eps}\in\cu$ s.t. $\eps e_{n} -a_{n}\geq -\eps e_{n}$, for
     all $n\in U_{\eps}$, which implies that, for all $\eps>0$,
     $\{n\in\bn : a_{n}<\eps e_{n}\}\in \cu$.  Hence, because $\cu$ is
     free, $U_{k} := \{n\geq k: a_{n} < \frac{1}{k} e_{n}\}\in\cu$. 
     Clearly $U_{k+1}\subset U_{k}$, $k\in\bn$, and $\cap_{k\in\bn}
     U_{k} =\emptyset$.  Set $G_{0}:= \bn\setminus U_{1}$, $G_{k}:=
     U_{k}\setminus U_{k+1}$, $k\in\bn$, and 
     $$
     b_{n}:= 
     \begin{cases}
	 \|a_{n}\|e_{n}& n\in G_{0} \\
	 \frac{1}{k}e_{n}& n\in G_{k}.
     \end{cases}
     $$
     This implies $\lim_{\cu} \|b_{n}\| = 0$, and $a_{n}-b_{n}\leq 0$, 
     $n\in\bn$, that is $a_{\cu} = \lim_{\cu} a_{n} = \lim_{\cu} 
     (a_{n}-b_{n}) \leq 0$, which is the thesis.
 \end{proof}
 
 \begin{Prop}\label{prop:1.16}
     Let $\{(X_{n},e_{n},L_{n})\}$ be a uniform sequence of Lip-normed
     order-unit spaces, $\cu$ a free ultrafilter.  Then the
     Lip-ultraproduct $(\ell^{\infty}_{L}(X,\cu),e_{\cu})$ is a
     Lip-normed order-unit space.
 \end{Prop}
 \begin{proof}
    It follows from Theorem \ref{thm:1.12} and Proposition
    \ref{restrUltra} that $(X_{\cu},e_{\cu})$ is a \lsp, with $\|
    a_{\cu}\|_{L,\cu} := \inf_{\{y_{n}\}\equiv \{a_{n}\}}
    \sup_{n}\|y_{n}\|_{L}$.  Then $(X_{\cu},e_{\cu})$ is an order-unit
    space.  Indeed, $\ell^{\infty}_{L}(X_{n})$ is a closed subspace of
    $\ell^{\infty}(X_{n})$, containing $e:=\{e_{n}\in
    X_{n}\}_{n\in\bn}$.  So $\ell^{\infty}_{L}(X_{n})\cap\cai_{\cu}$
    is a positively generated order-ideal of
    $\ell^{\infty}_{L}(X_{n})$, and arguing as in the previous
    Proposition, $X_{\cu}=\pi_{\cu}(\ell^{\infty}_{L}(X_{n})) =
    \ell^{\infty}_{L}(X_{n})/\ell^{\infty}_{L}(X_{n})\cap\cai_{\cu}$
    is Archimedean, therefore an order-unit space.
    
    Let us set $L(a_{\cu}) := \inf_{\l\in\br} \| a_{\cu}-\l
    e_{\cu}\|_{L,\cu}$. Then it follows from Proposition \ref{Prop1.1a} 
    that $L$ is a lower semicontinuous Lipschitz seminorm.
    Finally we prove that $\rho_{L}$ induces on $S(X_{\cu})$ the 
    weak$^{*}$-topology. Indeed, for $\o_{1},\o_{2}\in S(X_{\cu})$, we have
    \begin{align}\label{eq:rho=dualLnorm}
	\rho_{L}(\o_{1},\o_{2}) & = \sup_{a} 
	\frac{|\o_{1}(a)-\o_{2}(a)|}{L(a)} 
	= \sup_{a} \frac{|\o_{1}(a)-\o_{2}(a)|}{\inf_{\l} \|
	a_{\cu}-\l e_{\cu}\|_{L,\cu}}  \notag\\
	& = \sup_{a, \l} \frac{|\o_{1}(a-\l e_{\cu})-\o_{2}(a-\l
	e_{\cu})|}{\| a_{\cu}-\l e_{\cu}\|_{L,\cu}}  
	= \sup_{a} \frac{|\o_{1}(a)-\o_{2}(a)|}{\| a_{\cu}\|_{L,\cu}} \notag \\
	& = \|\o_{1}-\o_{2}\|'_{L}.
    \end{align}
    Therefore $\rho_{L}$ induces on $S(X_{\cu})$ the 
    weak$^{*}$-topology by Proposition \ref{prop:1.2}, and $L$ is a 
    Lip-seminorm.
 \end{proof}

 The seminorm $L$ in the previous Proposition can be obtained more 
 directly in terms of the seminorms $L_{n}$, as the following 
 Proposition shows.
 
 \begin{Prop}
     Let $\{(X_{n},e_{n},L_{n})\}$ and $\cu$ be as in the previous
     Proposition.  Then 
     \item{$(i)$} The Lip-seminorm on the Lip-ultraproduct of
     order-unit spaces gives back the Lip-norm on the Lip-ultraproduct
     of Lip-spaces, namely, for any $x_{\cu}$ in the ultraproduct,
         \begin{equation}\label{wayback}
	\|x_{\cu}\|_{L,\cu} = \max \{
	\frac{\|x_{\cu}\|}{R_{\cu}} , L_{\cu}(x_{\cu}) \}.
    \end{equation}    
     \item{$(ii)$} The Lip-seminorm on the Lip-ultraproduct is the
     quotient seminorm, namely
     \begin{equation}
	 L(x_{\cu}) = \inf_{[x_{n}]= x_{\cu}} \sup_{n}
	 L_{n}(x_{n}).
     \end{equation}
 \end{Prop}
 \begin{proof}
    Let us first observe that
    $$
    \lim_{\cu} R_{n} = R_{\cu} 
    = \sup \frac{\|x_{\cu}\|}{\|x_{\cu}\|_{L}},
    $$
    where we used Propositions \ref{prop:1.13}$(vi)$ and \ref{restrUltra}.     
    
    Let us now set $L_{\cu}(x_{\cu}) := \inf_{[x_{n}]= x_{\cu}}
    \sup_{n} L_{n}(x_{n})$.  We want to prove that, $\forall x_{\cu}\in
    X_{\cu}$, $\exists \{\tilde{x}_{n}\}\in\ell^{\infty}_{L}(X_{n})$
    s.t. $[\tilde{x}_{n}]=x_{\cu}$ and
    \begin{equation}\label{rappresentante}
	\lim_{\cu}\|\tilde{x}_{n}\|_{L} = \|x_{\cu}\|_{L,\cu} \qquad \lim_{\cu} 
	L_{n}(\tilde{x}_{n}) = L_{\cu}(x_{\cu}).
    \end{equation}
    
    Let $x_{\cu}\in X_{\cu}$, and, for any $k\in\bn$, choose sequences
    $x^{k}_{n}$ realising it and such that
    \begin{equation}\label{eq:c}
        L_{n}(x^{k}_{n})\leq (1+\frac1k)L_{\cu}(x_{\cu}), \quad n\in\bn
    \end{equation}
    As $\lim_{\cu} \frac{\|x^{k}_{n}\|}{R_{n}} =
    \frac{\|x\|_{\cu}}{R_{\cu}}$, there is $U_{k}\in\cu$ s.t.
    $\frac{\|x^{k}_{n}\|}{R_{n}} \leq (1+\frac1k)
    \frac{\|x\|_{\cu}}{R_{\cu}}$, $n\in U_{k}$.  Setting, if necessary,
    $$
    \widetilde{x^{k}_{n}} := 
    \begin{cases}
	x^{k}_{n}& n\in U_{k}\\
	0&n\not\in U_{k},
    \end{cases}
    $$
    we obtain 
    \begin{equation}\label{eq:a}
	\frac{\|\widetilde{x^{k}_{n}}\|}{R_{n}} \leq (1+\frac1k) 
    \frac{\|x\|_{\cu}}{R_{\cu}}, n\in\bn,
    \end{equation}
    and $\{\widetilde{x^{k}_{n}}\} \equiv \{x^{k}_{n}\}$, for all 
    $k\in\bn$. Therefore we can assume that $\{x^{k}_{n}\}$ have been 
    chosen in such a way that (\ref{eq:c}), (\ref{eq:a}) are satisfied.
    
    Using (\ref{quotientLip}) and $L_{n}(y_{n}) \leq \|y_{n}\|_{L}$, 
    we obtain $L_{\cu}(x_{\cu}) \leq \|x_{\cu}\|_{L}$.
    
    Set, for all $k\in\bn$, $V_{k} = \{ n \geq k : \|x^{i}_{n} -
    x^{j}_{n} \| \leq \frac1i, i\leq j\leq k\}$, $V_{0}:= 
    \bn\setminus V_{1}$, and then
    $\tilde{x}_{n}=\frac{k}{k+1} x^{k}_{n}$, $n\in V_{k}\setminus
    V_{k+1}$.  Then, $[\tilde{x}_{n}]=x_{\cu}$.  Moreover, for
    $k,n\in\bn$, we have, using (\ref{eq:c}), (\ref{eq:a}),
    \begin{align*}
	\|x^{k}_{n}\|_{L} & = \max \left\{ L_{n}(x^{k}_{n}), 
	\frac{\|x^{k}_{n}\|}{R_{n}} \right\} \\
	& \leq (1+\frac1k) \max \left\{ L_{\cu}(x_{\cu}), 
	\frac{\|x_{\cu}\|_{\cu}}{R_{\cu}} \right\} \leq (1+\frac1k) 
	\|x_{\cu}\|_{L},
    \end{align*}
    so that, for $k\in\bn$, $\ell\in\bn$, $\ell\geq k$, $n\in 
    V_{\ell}\setminus V_{\ell+1}$, we get $\| \tilde{x}_{n} \|_{L} = 
    \frac{\ell}{\ell+1} \|x^{\ell}_{n}\|_{L} \leq \|x_{\cu}\|_{L}$, 
    which implies $\| \tilde{x}_{n} \|_{L} \leq \|x_{\cu}\|_{L}$, for 
    $n\in V_{k}$, and $\lim_{\cu} \| \tilde{x}_{n} \|_{L} \leq \|x_{\cu}\|_{L}$.
    As the opposite inequality is always true, we obtain
    $$
	\lim_{\cu} \| \tilde{x}_{n} \|_{L} = \|x_{\cu}\|_{L}.
    $$
    Finally, from (\ref{eq:c}), for $k\in\bn$, $\ell\in\bn$, $\ell\geq k$, $n\in 
    V_{\ell}\setminus V_{\ell+1}$, we get $L_{n}(\tilde{x}_{n}) = 
    \frac{\ell}{\ell+1} L_{n}(x^{\ell}_{n}) \leq L_{\cu}(x_{\cu})$, 
    which implies $L_{n}(\tilde{x}_{n})  \leq L_{\cu}(x_{\cu})$, $n\in 
    V_{k}$, and $\lim_{\cu} L_{n}(\tilde{x}_{n})$  $\leq L_{\cu}(x_{\cu})$. 
    As the opposite inequality is always true, we obtain
    $$
	\lim_{\cu} L_{n}(\tilde{x}_{n}) = L_{\cu}(x_{\cu}),
    $$
    and we have proved (\ref{rappresentante}).
    
    As a consequence, we get equation (\ref{wayback}):
    \begin{align*}
	\|x_{\cu}\|_{L,\cu} & = \lim_{\cu} \|\wt{x}_{n}\|_{L} =
	\lim_{\cu} \max \{ \frac{\|\wt{x}_{n}\|}{R_{n}} ,
	L_{n}(\wt{x}_{n}) \} \notag\\
	& = \max \{\lim_{\cu} \frac{\|\wt{x}_{n}\|}{R_{n}} ,
	\lim_{\cu} L_{n}(\wt{x}_{n}) \} = \max \{
	\frac{\|x_{\cu}\|}{R_{\cu}} , L_{\cu}(x_{\cu}) \}.
    \end{align*}    
    Let us now denote by $\l_{n}$ the constant for which
    $\|\tilde{x}_{n}\|_{0}=\|\tilde{x}-\l_{n}e_{n}\|$.  Since
    $\{\tilde{x}_{n}\}$ is norm bounded, $\{\l_{n}\}$ is bounded,
    hence $\lim_{\cu}\l_{n}=\l_{\cu}\in\br$.  Then
    \begin{align*}
	\|x_{\cu}\|_{0} & = \inf_{\l} \|x_{\cu} -\l e_{\cu}\| \leq
	\|x_{\cu} -\l_{\cu} e_{\cu}\| 
	= \lim_{\cu}  \| \wt{x}_{n} -\l_{n} e_{n} \|\\
	& = \lim_{\cu} \|\wt{x}_{n}\|_{0} \leq \lim_{\cu} R_{n}
	L_{n}(\wt{x}_{n}) = R_{\cu} L_{\cu}(x_{\cu}).
    \end{align*}
    Therefore, using (\ref{wayback}) for the vector $x_{\cu} -\l 
    e_{\cu}$ and the inequality above,
    \begin{align*}
	L(x_{\cu}) & = \inf_{\l} \|x_{\cu} -\l e_{\cu} \|_{L,\cu} =
	\max \{ \frac{\inf_{\l} \|x_{\cu} -\l e_{\cu}
	\|_{\cu}}{R_{\cu}} , L_{\cu}(x_{\cu}) \} \\
	& = \max \{ \frac{\|x_{\cu}\|_{0}}{R_{\cu}} , L_{\cu}(x_{\cu})
	\} = L_{\cu}(x_{\cu}),
    \end{align*}
    concluding the proof.
\end{proof}
 
 Now we can prove the analogue of Theorem \ref{thm:dualUltraprod}.

\begin{Thm}\label{Thm:UPOUS}
    Given a uniform sequence $\{(X_{n},e_{n},L_{n})\}$ of Lip-normed
    order-unit spaces, the ultraproduct $\ell^{\infty}(X_{n}',\cu)$ of
    the dual spaces projects on the dual $(\ell^{\infty}_{L} (X_{n},
    \cu))'$ of the Lip-ultraproduct.  Moreover, any state on
    $\ell^{\infty}_{L} (X_{n}, \cu)$ can be represented by an element
    of $\ell^{\infty}(X_{n}',\cu)$ given by sequences of states.
\end{Thm}
\begin{proof}
    Only the last part needs a proof, which is similar to that of
    Theorem \ref{thm:dualUltraprod}, so that we only indicate the 
    small difference.
    
    Given $\eps>0$, let us choose the subspaces $V_{n}\subset X_{n}$
    as in the proof of the cited Theorem, but with the further request 
    that $e_{n}\in V_{n}$, for any $n\in\bn$.
    
    Now take any $\f \in S(X_{\cu})$, follow the proof of the cited
    Theorem until you get elements $\f'_{n}\in X'_{n}$, and set $\f'
    := [\f'_{n}]\in \ell^{\infty}(X'_{n},\cu)$.  In this case, since
    $V_{n}\ni e_{n}$, $\|\f'_{n}\|=1$, and recall that, by
    construction, $\f$ and $\pi(\f')$ coincide on $V_{\cu}$. 
    Therefore $\lim_{\cu}\f'_{n}(e_{n})=\pi(\f')(e_{\cu})=1$.  Now we
    may decompose $\f'_{n}=\a_{n}\psi^{1}_{n}-\b_{n}\psi^{2}_{n}$
    where $\a_{n}+\b_{n}=1$, $\a_{n}\geq0$, $\b_{n}\geq0$, and
    $\psi^{i}_{n}$ are states (\cite{Alfsen}, II.1.14).  Therefore we
    obtain $\a_{n}\to1$ and $\b_{n}\to0$, implying that
    $[\psi^{1}_{n}]=\f'$, namely can be realised via sequences of
    states.  The proof continues as in the cited Theorem.
\end{proof}

Let us recall that \cite{BHbook}, given a sequence $(X_{n},d_{n})$ of metric spaces
with uniformly bounded radius, and $\cu$ a free ultrafilter on $\bn$,
the ultralimit $(X_{\cu},d_{\cu})$ is defined as the space of
equivalence classes $[x_{n}]$, $x_{n}\in X_{n}$, with distance
$d_{\cu}([x_{n}],[x'_{n}]) = \lim_{\cu}d(x_{n},x'_{n})$, and it
follows that $[x_{n}]=[x'_{n}]$ when they have zero distance. 
According to Proposition \ref{Prop1.6}, we have $\rho_{L_{\cu}}
(\f_{\cu}, \psi_{\cu}) = \lim_{\cu} \rho_{L} (\f_{n}, \psi_{n})$,
therefore we get the following.

\begin{Cor}
    Let $\{(X_{n},e_{n},L_{n})\}$ be a uniform sequence of Lip-normed
    order-unit spaces.  The state space of the Lip-ultraproduct can be
    isometrically identified with the ultralimit of the approximating
    state spaces.
\end{Cor}

 Let us now recall Rieffel's notion of quantum Gromov-Hausdorff 
 convergence \cite{Rieffel00}. 
 
 Let $(X,e_{X},L_{X})$, $(Y,e_{Y},L_{Y})$ be Lip-normed order-unit
 spaces.  Denote by ${\mathcal M}(L_X,L_Y)$ the set of lower
 semicontinuous Lip-seminorms on $X \oplus Y$ which induce $L_X$ and
 $L_Y$ on $X$, $Y$ respectively.  Any $L \in {\mathcal M}(L_X,L_Y)$
 gives rise to a metric $\rho_L$, on $S(X \oplus Y)$.  Therefore,
 identifying $S(X)$ and $S(Y)$ with (closed, convex) subsets of $S(X
 \oplus Y)$, we can consider the Hausdorff distance between them
 w.r.t. $\rho_L$, namely $\rho^{H}_L(S(X),S(Y))$.  We define the {\em
 quantum Gromov--Hausdorff distance} between $X$ and $Y$ by
 \begin{equation}\label{qGHdist}
 \dist(X,Y) = \inf\{\rho^{H}_L(S(X),S(Y)): L \in {\mathcal
 M}(L_X,L_Y)\}.
 \end{equation}

 \begin{Thm}\label{thm:orderunitIsComplete}
     Let $\{(X_{n},e_{n},L_{n})\}$ be a Cauchy sequence of Lip-normed
     order-unit spaces.  Then, for any free ultrafilter $\cu$, the
     Lip-normed Lip-ultraproduct
     $(\ell^{\infty}_{L}(X_{n},\cu),e_{\cu},L_{\cu})$ is the limit of
     the sequence.
 \end{Thm}
 \begin{proof}
     Let $\eps>0$ be given, and let $n_{\eps}\in\bn$ be s.t. for all
     $m,n>n_{\eps}$ there is $L_{mn}\in \cam(X_{m}, X_{n})$ s.t.
     $\r^{H}_{L_{mn}} (S(X_{n}), S(X_{m})) <\eps$.  Observe that,
     having fixed $n>n_{\eps}$, the Lip-ultraproduct of the spaces
     $\{X_{n}\oplus X_{i}\}_{i\in\bn}$ naturally identifies with
     $X_{n}\oplus X_{\cu}$.  Therefore, $X_{n}\oplus X_{\cu}$ inherits
     a Lip-seminorm $L_{n\cu}$ with respect to which $S(X_{\cu})
     \subset \ov{B}_{\eps}(S(X_{n}))$ and $S(X_{n}) \subset
     \ov{B}_{\eps}(S(X_{\cu}))$.
     
     Indeed, if $\o\in S(X_{n})$, then, for all $m>n_{\eps}$, there is
     $\f_{m}\in S(X_{m})$ s.t. $\rho_{L_{mn}}(\o,\f_{m})<\eps$.  Set
     $\f_{\cu}(x_{\cu}):= \lim_{\cu}\f_{m}(x_{m})$, $[x_{m}]=x_{\cu}$
     (see Proposition \ref{Prop1.6}) so that $\f_{\cu}\in S(X_{\cu})$
     and, by (\ref{eq:rho=dualLnorm}) and Proposition \ref{Prop1.6},
     $$
     \rho_{L_{n\cu}}(\o,\f_{\cu}) = \lim_{m\to\cu}
     \rho_{L_{mn}}(\o,\f_{m}) \leq \eps.
     $$
     
     Viceversa, let $\f\in S(X_{\cu})$, and choose, by Theorem
     \ref{Thm:UPOUS}, $\f_{m}\in S(X_{m})$, s.t. $\f_{\cu}(x_{\cu}):=
     \lim_{\cu}\f_{m}(x_{m})$, $[x_{m}]=x_{\cu}$, and let, for
     $m>n_{\eps}$, $\o_{m}\in S(X_{n})$ be s.t.
     $\r_{L_{mn}}(\f_{m},\o_{m}) <\eps$.  Set $\o:=\lim_{m\to\cu}
     \o_{m}\in S(X_{n})$.  Then, $\r_{L_{n\cu}}(\o,\f) =
     \lim_{m\to\cu} \r_{L_{mn}}(\o,\f_{m}) \leq \lim_{m\to\cu}
     \r_{L_{mn}}(\o_{m},\f_{m}) + \r_{L_{n}}(\o_{m},\o) \leq\eps$.
 \end{proof}

\section{Operator Systems}\label{sec:3}
 
We begin by describing our operator system framework.  For references
see \cite{Paulsen}.  

\begin{Dfn}\label{def:3.1}
    An operator system $X$ is a complex vector space with a
    conjugate linear involution $^{*}:x\in X\to x^{*}\in X$,
    satisfying
    \itm{i} $X$ is matrix ordered, $i.e.$ 
    \itm{i'} for any $p\in\bn$, there is a proper cone
    $M_{p}(X)_{+}\subset M_{p}(X)_{h}$, where the subscript $h$ refers 
    to hermitian elements
    \itm{i''} for any $p,q\in\bn$, $A\in M_{qp}(\bc)$, $A^{*}
    M_{q}(X)_{+} A \subset M_{p}(X)_{+}$
    \itm{ii} $X$ has a matrix order-unit, $i.e.$ there is $e\in X_{h}$
    s.t., with $e^{p} := \diag(e,\ldots,e)\in M_{p}(X)_{+}$, for any
    $x\in M_{p}(X)_{h}$, there is $r>0$ s.t. $x+re^{p}\in M_{p}(X)_{+}$
    \itm{iii} the matrix order-unit $e$ is Archimedean, $i.e.$ if
    $x\in M_{p}(X)$ is s.t. $x+re^{p}\in M_{p}(X)_{+}$, for all $r>0$,
    then $x\in M_{p}(X)_{+}$.
\end{Dfn}

 Given operator systems $X$ and $Y$ we say that a linear map $\varphi
 : X \to Y$ is {\em $n$-positive} if the map $\id_n \otimes\varphi :
 M_n \otimes X \to M_n \otimes Y$ is positive, and if $\id_n
 \otimes\varphi$ is positive for all $n\in\mathbb{N}$ then we say that
 $\varphi$ is {\em completely positive}.  A completely positive
 (resp.\ unital completely positive) linear map will be referred to as
 a {\em c.p.}\ (resp.\ {\em u.c.p.}) map.  If $\varphi : X \to Y$ is a
 unital $m$-positive map with $m$-positive inverse for $m=1, \dots ,
 n$ then $\varphi$ is a {\em unital $n$-order isomorphism}, and if
 $\varphi$ is u.c.p.\ with c.p.\ inverse then $\varphi$ is a {\em
 unital complete order isomorphism}.

 We denote by $UCP_n (X)$ the collection of all u.c.p.\ maps from $X$
 into $M_n$ (the {\em matrix state spaces}).

 Following Kerr \cite{Kerr}, we introduce Lip-norms and matricial
 distances on operator systems.  By a {\em Lip-normed operator system}
 we mean a pair $(X,L)$ where $X$ is a complete operator system and
 $L$ is a lower semicontinuous Lip-seminorm on $X$ satisfying
 $L(x^{*})=L(x)$.  If $X$ is a unital $C^*$-algebra then we will also
 refer to $(X,L)$ as a {\em Lip-normed unital $C^*$-algebra}.

\begin{Dfn}\label{D-ucpmetric}
    Let $(X,L)$ be a Lip-normed operator system and $p\in\mathbb{N}$. 
    We define the metric $\rho_{L}$ on $UCP_p (X)$ by 
    $$ 
    \rho_{L} (\varphi , \psi ) = \sup_{L(x)\leq 1} \|
    \varphi (x) - \psi (x) \|
    $$
    for all $\varphi , \psi\in UCP_p (X)$,
\end{Dfn}

 Let $(X,L_X )$ and $(Y,L_Y )$ be Lip-normed operator systems.  We
 denote by $\mathcal{M}(L_X , L_Y )$ the collection of lower
 semicontinuous Lip-seminorms on $X \oplus Y$ which induce $L_X$ and
 $L_Y$ via the quotient maps onto $X$ and $Y$, respectively.

 Let $L\in\mathcal{M}(L_X , L_Y )$.  Since the projection map $X\oplus
 Y\to X$ is u.c.p., by \cite{Kerr}, we obtain an isometry $UCP_p
 (X)\to UCP_p (X\oplus Y)$ with respect to $\rho_{L_X}$ and $\rho_L$. 
 Similarly, we also have an isometry $UCP_p (Y)\to UCP_p (X\oplus Y)$. 
 For notational simplicity we will thus identify $UCP_p (X)$ and
 $UCP_p (Y)$ with their respective images under these isometries.

\begin{Dfn}\label{D-dist}
    Let $(X,L_X )$ and $(Y,L_Y )$ be Lip-normed operator systems.  For
    each $p\in\mathbb{N}$ we define the $p$-distance 
    $$
    \dist_{p} (X,Y) = \inf_{L\in\mathcal{M}(L_X , L_Y )}
    \rho^{H}_{L} (UCP_p (X) , UCP_p (Y))
    $$
    where $\rho^{H}_{L}$ denotes Hausdorff distance with
    respect to the metric $\rho_{L}$.  We also define the 
    complete quantum Gromov-Hausdorff distance
    $$ 
    \dist_{\infty} (X,Y) = \inf_{L\in\mathcal{M}(L_X , L_Y )} \,\,
    \sup_{p\in\mathbb{N}}\, \rho^{H}_{L} (UCP_p (X) , UCP_p
    (Y)) .
    $$
\end{Dfn}
 
 \begin{Prop}\label{prop:2.6}
     Let $\{(X_{n},e_{n})\}$ be operator systems, $\cu$ a free ultrafilter. 
     Then the ultraproduct $(\ell^{\infty}(X_{n},\cu),e_{\cu})$ is an
     operator system.
 \end{Prop}
 \begin{proof}
     Denote $X_{\cu}:= \ell^{\infty}(X_{n},\cu)$. 
     It follows from Definition \ref{def:3.1} $(i'),(ii),(iii)$, that, 
     for any $p\in\bn$, $(M_{p}(X_{n}),e^{p}_{n})$ is a complete order-unit 
     space, so that $(M_{p}(X_{\cu})\equiv 
     \ell^{\infty}(M_{p}(X_{n}),\cu), e^{p}_{\cu})$ is a complete order-unit 
     space, by Proposition \ref{prop:1.15}. Finally, for any 
     $p,q\in\bn$, $A\in M_{qp}(\bc)$, from $A^{*}
    M_{q}(X_{n})_{+} A \subset M_{p}(X_{n})_{+}$ it follows $A^{*}
    M_{q}(X_{\cu})_{+} A \subset M_{p}(X_{\cu})_{+}$. Therefore 
    $(X_{\cu},e_{\cu})$ is a complete operator system.
 \end{proof}

 \begin{Prop}
     Let $\{(X_{n},e_{n},L_{n})\}$ be a uniform sequence of Lip-normed
     operator systems, $\cu$ a free ultrafilter.  Then the
     Lip-ultraproduct $(\ell^{\infty}_{L}(X_{n},\cu),e_{\cu},L_{\cu})$
     is a Lip-normed operator system.
 \end{Prop}
 \begin{proof}
     It follows from Propositions \ref{prop:2.6} and \ref{prop:1.16}. 
 \end{proof}

 \begin{Prop}\label{Prop6}
     Let $\{(X_{n},e_{n},L_{n})\}$ be a uniform sequence of Lip-normed
     operator systems, $\cu$ a
     free ultrafilter.  Let $p\in\bn$, $\{\s_{n}\in UCP_{p}(X_{n})\}$,
     $\{\t_{n}\in UCP_{p}(X_{n})\}$.  Define
     $\s_{\cu}(a_{\cu}) := \lim_{\cu} \s_{n}(a_{n})$, $a_{\cu}=[a_{n}]\in
     \ell^{\infty}_{L}(X_{n},\cu)$, and $\t_{\cu}$ analogously.  
     Then $\s_{\cu},\t_{\cu}$ are well defined and belong to 
     $UCP_{p}(X_{\cu})$, and
     $$
     \r_{L_{\cu}}(\s_{\cu},\t_{\cu}) = \lim_{\cu} \r_{L_{n}}(\s_{n},\t_{n}).
     $$
 \end{Prop}
 \begin{proof}
     The first part is as in Proposition \ref{Prop1.6}. Moreover
     \begin{align*}
	 \lim_{\cu} \r_{L_{n}}(\s_{n},\t_{n}) & = \lim_{\cu}
	 \sup_{x_{n} \in X_{n}}
	 \frac{\|\s_{n}(x_{n})-\t_{n}(x_{n})\|}{ L_{n}(x_{n}) } =
	 \sup_{x\in \ell^{\infty}_{L}(X_{n})} \lim_{\cu}
	 \frac{\|\s_{n}(x_{n})-\t_{n}(x_{n})\|}{ L_{n}(x_{n}) } \\
	 & = \sup_{x_{\cu}\in X_{\cu}}
	 \sup_{[x_{n}]=x_{\cu}}
	 \frac{\|\s_{\cu}(x_{\cu})-\t_{\cu}(x_{\cu})\|}{\lim_{\cu}
	 L_{n}(x_{n}) } \\
	 & = \sup_{x_{\cu}\in X_{\cu}}
	 \frac{\|\s_{\cu}(x_{\cu})-\t_{\cu}(x_{\cu})\|}{
	 L_{\cu}(x_{\cu}) } = \r_{L_{\cu}}(\s_{\cu},\t_{\cu}),
     \end{align*}
     where in the last but one equality we used
     (\ref{rappresentante}), and the consideration at the end of the
     proof of Proposition \ref{Prop1.6} applies.
 \end{proof}

 \begin{Thm}\label{thm:2.12}
     Let $\{(X_{n},e_{n},L_{n})\}$ be a Cauchy sequence of Lip-normed
     operator systems.  Then $(X_{\cu},e_{\cu},L_{\cu})$ is its limit,
     for any free ultrafilter $\cu$.
 \end{Thm}
 \begin{proof}
     It is similar to the proof of Theorem \ref{thm:orderunitIsComplete}, 
     by making use of Proposition \ref{Prop6}, and the analogue of 
     Theorem \ref{Thm:UPOUS}.
 \end{proof}

 \section{C$^*$-algebras}

 \subsection{The problem of completeness}
 
 Let us consider the space of equivalence classes of Lip-normed unital
 C$^{*}$-algebras, endowed with one of the pseudo-distances
 $\dist_{p}$, $p\in\bn\cup\{\infty\}$.  Kerr showed \cite{Kerr} that
 for $p\geq2$ it is indeed a distance, namely that if
 $\dist_{p}(\ca,\cb)=0$ then $\ca$ and $\cb$ are Lip-isometric
 $^{*}$-isomorphic C$^{*}$-algebras.
  
 Our aim is to study the completeness of the equivalence classes of
 C$^{*}$-algebras endowed with the metrics $\dist_{p}$.  When
 $\dist_{\infty}$ is considered, the limit of a Cauchy sequence exists
 as an operator system.  The result of Kerr implies that, on such a
 space, the C$^{*}$-structure, i.e. a product w.r.t. which the norm is
 a C$^{*}$-norm, is unique, if it exists.
 However, besides the mere question of existence of such a product, we
 are interested in products which are approximated by the products of
 the approximating algebras. 

 A first attempt in this respect has been made by David Kerr and
 Hanfeng Li \cite{Kerr,Li2}, who introduced the concept of $f$-Leibniz
 property, showing that if all algebras in a Cauchy sequence enjoy the
 $f$-Leibniz property for the same function $f$, then the limit space
 inherits a product structure (satisfying the $f$-Leibniz property).
 
 We observe however that realising the limit space as a
 Lip-ultraproduct allows a much more stringent characterization of the
 cases in which the product structure is inherited by the limit space.

 Indeed, when realising the limit as a Lip-ultraproduct, one would 
 like to set
 \begin{equation}\label{inherprod}
     [x_{n}]\ [y_{n}]=[x_{n}y_{n}].
 \end{equation}
 Unfortunately it is not true in general that $[x_{n}y_{n}]$ belongs
 to the Lip-ultraproduct, namely has finite Lip-norm or at least can
 be approximated in norm by elements with finite Lip-norm.  In other
 words, while (\ref{inherprod}) defines a product on
 $\ell^{\infty}(\ca_{n},\cu)$, it is not always true that
 $\ell^{\infty}_{L}(\ca_{n},\cu)$ is a subalgebra of
 $\ell^{\infty}(\ca_{n},\cu)$.

 We then introduce the following
 
 \begin{Dfn}
     Let $\{\ca_{n}\}_{n\in\bn}$ be a Cauchy sequence of Lip-normed
     unital C$^{*}$-algebras w.r.t. the $\dist_{p}$ metrics.  If $\cu$
     is a free ultrafilter on $\bn$, we say that the Lip-ultraproduct
     $\ell^{\infty}_{L}(\ca_{n},\cu)$ inherits the C$^{*}$-structure
     if it is a sub-algebra of $\ell^{\infty}(\ca_{n},\cu)$.  In
     general, we say that the limit inherits the C$^{*}$-structure if
     $\ell^{\infty}_{L}(\ca_{n},\cu)$ does, for some free ultrafilter
     $\cu$.
 \end{Dfn}
 
 \begin{Prop}\label{C*ultra}
     Let $\{\ca_{n}\}_{n\in\bn}$ be a Cauchy sequence of Lip-normed
     unital C$^{*}$-algebras w.r.t. the $\dist_{p}$ metrics, and
     suppose $\ell^{\infty}_{L}(\ca_{n},\cu)$ inherits the
     C$^{*}$-structure for a suitable free ultrafilter $\cu$.  Then
     the sequence $\{\ca_{n}\}_{n\in\bn}$ converges to the \cst{}
     $\ell^{\infty}_{L}(\ca_{n},\cu)$.
 \end{Prop}
 
 \begin{proof} 
     Cf. the proofs of Theorems \ref{thm:orderunitIsComplete}, \ref{thm:2.12}.
 \end{proof}
 
 As we shall see in Subsection \ref{examples}, the general situation
 is as ugly as possible: there are Cauchy sequences for which the
 limit is not a C$^{*}$-algebra, and even Cauchy sequences for which
 the limit can be endowed with a C$^{*}$-product, but this is not
 inherited from the approximating C$^{*}$-algebras.
 
 \begin{Thm}\label{noncompl}
     \item{$(i)$} The space of equivalence classes of Lip-normed
     unital C$^{*}$-algebras, endowed with the distance $\dist_{p}$,
     $p\geq2$, is not complete.  \item{$(ii)$} There exist sequences
     $(\ca_{n},L_{n})$ converging to a Lip-normed unital
     {C$^{*}$-algebra} $(\ca,L)$ for which the C$^{*}$-structure is
     not inherited.
 \end{Thm}

 We are not able to characterise the Cauchy sequences for which the limit
 admits a C$^{*}$-product, but we can characterise those for which the 
 C$^{*}$-product is inherited. Our condition may be seen as a 
 generalisation of the Kerr-Li condition.

 \begin{Dfn}
     We say that the pair $(\ca,L)$ consisting of a unital
     {C$^{*}$-algebra} and a seminorm is a Lip-normed unital
     {C$^{*}$-algebra} if $L$ is a lower semicontinuous Lip-seminorm
     according to Section \ref{sec:3}.  $(\ca,L)$ will be called {\it
     quasi Lip-normed} if we drop the assumption that Lip-elements are
     dense, but assume that they generate $\ca$ as a C$^{*}$-algebra.
 \end{Dfn}

 Given a quasi Lip-normed unital {C$^{*}$-algebra}, we consider the
 function
 $$
 \eps(r)=\sup_{\nl{x}\leq 1}\inf_{\nl{y}\leq r}\|y-x^{*}x\|,
 $$
 where $\nl{\cdot}$ denotes the Lipschitz norm defined in subsection 
 \ref{subsec:1.3}.
 
 \begin{Lemma}
     The quasi Lip-normed unital {C$^{*}$-algebra} $(\ca,L)$ is
     Lip-normed if and only if
     \begin{equation}\label{epsto0}
	 \lim_{r\to\infty}\eps(r)=0.
     \end{equation}
 \end{Lemma}
 
 \begin{proof}
     Assume Lip-elements are dense. This means that, for any $\eps>0$, 
     the open sets 
     $$
     \Q(\eps,r)=\bigcup_{\|x\|_{L}\leq r}B(x,\eps),\quad r>0
     $$
     give an open cover of $\ca$. Since $\{x^{*}x:\|x\|_{L}\leq1\}$ 
     is compact, we may extract a finite subcover, hence 
     $\forall\eps>0$ $\exists r>0$ s.t. $\{x^{*}x:\|x\|_{L}\leq1\}\subset 
     \Q(\eps,r)$, or, equivalently, $\forall\eps>0$ $\exists 
     r>0$ s.t. $\eps(r)<\eps$, proving one implication.
     
     Now assume $\eps(r)\to0$.  This implies that for any Lip-element
     $x$, $x^{*}x$ can be arbitrarily well approximated (in norm) by
     Lip-elements.  Since $xy$ can be written as a linear combination
     of $x^{*}x$, $y^{*}y$, $(x+y)^{*}(x+y)$ and $(x+iy)^{*}(x+iy)$,
     we may conclude that products of Lip-elements can be arbitrarily
     well approximated (in norm) by Lip-elements.  Now take two
     norm-one elements $x$ and $y$ in the norm closure of the space of
     Lip-elements.  Choose two Lip-elements $x_{\eps}$, $y_{\eps}$,
     still with norm one, such that $\|x-x_{\eps}\|<\eps$,
     $\|y-y_{\eps}\|<\eps$, and then a Lip-element $z$ such that
     $\|x_{\eps}y_{\eps}-z\|<\eps$.  We get
     $$
     \|xy-z\|\leq\|xy-x_{\eps}y_{\eps}\|+\|x_{\eps}y_{\eps}-z\|\leq3\eps.
     $$
     This means that the norm closure of the space of Lip-elements is
     an algebra, hence a C$^{*}$-algebra.  By definition of quasi
     Lip-normed unital C$^{*}$-algebra, such closure coincides with
     $\ca$.
 \end{proof}
 
 Let us now compare condition (\ref{epsto0}) with the $f$-Leibniz 
 condition. Let us recall that $(\ca,L)$ satisfies the $f$-Leibniz 
 condition w.r.t. a given continuous 4-variable function $f$ if
 $$
 L(ab)\leq f(L(a),L(b),\|a\|,\|b\|),\qquad a,b\in\ca.
 $$
 
 \begin{Prop}\label{f-Leibniz}
     Let $(\ca,L)$ be a quasi Lip-normed unital {C$^{*}$-algebra}. 
     The following are equivalent: \item{$(i)$} $(\ca,L)$ satisfies
     the $f$-Leibniz condition w.r.t. some function $f$ \item{$(ii)$}
     $(\ca,L)$ satisfies the condition
     $$
     \|ab\|_{L}\leq C\|a\|_{L}\|b\|_{L},\qquad a,b\in\ca
     $$
     for some constant $C$
     \item{$(iii)$} the function $\eps(r)$ defined above is 
     zero for $r$ large enough.
 \end{Prop}
 
 \begin{proof}
     Clearly $(ii)\imply (i)$, since $\|a\|_{L} =
     \max\{R^{-1}\|a\|,L(a)\}$, with $R$ the radius of the state space. 
     Conversely, if we set
     $$
     K=\sup_{\|a\|_{L}\leq1,\|b\|_{L}\leq1}f(L(a),L(b),\|a\|,\|b\|),
     $$
     and observe that $K$ is finite by compactness, we get
     $$
     \|ab\|_{L}= \max\{R^{-1}\|ab\|,L(ab)\}\leq \max\{R,K\}\|a\|_{L}\|b\|_{L}.
     $$
     
     Now let us observe that $(iii)$ means that $\eps(r_{0})=0$ for
     some $r_{0}$, namely $\sup_{\|x\|_{L}\leq1}\|x^{*}x\|_{L}\leq
     r_{0}$ or, equivalently, $\|x^{*}x\|_{L}\leq r_{0}\|x\|_{L}^{2}$
     for any $x$.  The latter is clearly equivalent to property
     $(ii)$.
 \end{proof}
 
 Now we characterise the existence of an inherited C$^{*}$-structure. 
 Indeed, giving a uniform sequence $\ca_{n}$ of C$^{*}$-algebras with
 Lip-norms and a free ultrafilter $\cu$, we can construct the
 inclusions $\ell^{\infty}_{L}(\ca_{n},\cu) \subset \cb_{\cu} \subset
 \ell^{\infty}(\ca_{n},\cu)$, where $\cb_{\cu}$ denotes the
 C$^{*}$-algebra generated by $\ell^{\infty}_{L}(\ca_{n},\cu)$.  By
 the properties proved above, $\cb_{\cu}$ is a quasi Lip-normed unital
 C$^{*}$-algebra.

 \begin{Prop}
     Let $\{(\ca_{n},L_{n})\}_{n\in\bn}$ be a Cauchy sequence of
     Lip-normed unital C$^{*}$-algebras, with functions $\eps_{n}$,
     and let $\cb_{\cu}$ the quasi Lip-normed unital {C$^{*}$-algebra}
     defined above, with function $\eps_{\cu}$.  Then
     $$
     \eps_{\cu}(r)=\lim_{\cu}\eps_{n}(r).
     $$
 \end{Prop}
 
 \begin{proof}
     Given $r>0$, $n\in\bn$, let $x_{n},y_{n}\in\ca_{n}$ realise the
     worst element with Lip-norm $\leq1$ and the best approximation of
     $x^{*}_{n}x_{n}$ with Lip-norm $\leq r$ respectively, hence
     $\|x_{n}^{*}x_{n}-y_{n}\|=\eps_{n}(r)$, and then set $x =
     \lim_{\cu}x_{n} , y = \lim_{\cu}y_{n} , \eps(r) =
     \lim_{\cu}\eps_{n}(r)$.  This implies that
     $\|x^{*}x-y\|=\eps(r)$.  An element $\tilde{y} \in
     \ell^{\infty}_{L}(\ca_{n},\cu)$, $\|\tilde{y}\|_{L}\leq r$, giving
     the best approximation of $x^{*}x$, could be obtained as
     $\tilde{y} = \lim_{\cu}\tilde{y}_{n}$, with
     $\|\tilde{y}_{n}\|_{L}\leq r$, as shown in the proof of Lemma
     \ref{Lem:1.9}.  Since $\eps_{n}(r) \leq
     \|x_{n}^{*}x_{n}-\tilde{y}_{n}\| \to_{\cu}
     \|x^{*}x-\tilde{y}\|\leq \eps_{\cu}(r)$, we get
     $\eps(r)\leq\eps_{\cu}(r)$.
          
     Conversely, let $x\in\ell^{\infty}_{L}(\ca_{n},\cu)$ realise the
     worst element with Lip-norm $\leq1$, and, as above, obtain it as
     $x=\lim_{\cu}x_{n}$, $\|x_{n}\|_{L}\leq 1$.  Then let $y_{n}$ be the
     best approximation of $x^{*}_{n}x_{n}$ with Lip-norm $\leq r$ , hence
     $\|x^{*}_{n}x_{n}-y_{n}\|\leq\eps_{n}(r)$.  Setting
     $y=\lim_{\cu}y_{n}$, we get $\|y\|_{L}\leq r$ and
     $\eps_{\cu}(r)\leq\|x^{*}x-y\|\leq\eps(r)$.
 \end{proof}
 
 \begin{Cor}\label{Cor:inherited}
     Let $\{(\ca_{n},L_{n})\}_{n\in\bn}$ be a Cauchy sequence of
     Lip-normed unital C$^{*}$-algebras, with functions $\eps_{n}$. 
     The following are equivalent: \item{$(i)$} the limit inherits a
     C$^{*}$-structure
     \item{$(ii)$} $\displaystyle{ \lim_{r\to\infty}}
     \displaystyle{\lim_{\cu}\eps_{n}(r) = 0}$ for some free ultrafilter
     $\cu$ 
     \item{$(iii)$} there exists a subsequence $n_{k}$ such
     that
     $$
     \lim_{r\to\infty} \limsup_{k}\ \eps_{n_{k}}(r) = 0.
     $$
 \end{Cor}
 
 \begin{proof}
     By the results above, $(ii)$ amounts to saying that the quasi
     Lip-normed unital C$^{*}$-algebra{ } $\cb_{\cu}$ is indeed
     Lip-normed, hence coincides with
     $\ell^{\infty}_{L}(\ca_{n},\cu)$, which is therefore a
     C$^{*}$-algebra.
     
     $(iii)\implies (ii)$ For any free ultrafilter $\cu$ such that
     $\{n_{k}:k\in\bn\}\in\cu$, we have $\displaystyle{
     \lim_{r\to\infty}} \displaystyle{\lim_{\cu}\eps_{n}(r) = 0}$.
     
     $(ii)\implies (iii)$ Choose a sequence $\{ n^{1}_{k} \}_{ k\in\bn
     } \in \cu$ such that $\exists\lim_{k}\eps_{n^{1}_{k}}(1) =
     \eps_{\cu}(1)$, and then, inductively, $\{n^{j}_{k}\}_{k\in\bn}
     \in \cu$ as a subsequence of $n^{j-1}_{k}$ such that $\exists
     \lim_{k}\eps_{n^{j}_{k}}(j) = \eps_{\cu}(j)$.  For the diagonal
     subsequence $n_{k}:=n^{k}_{k}$, we get $\lim_{k}\eps_{n_{k}}(j) =
     \eps_{\cu}(j)$ for any $j$.  Then
     $$\limsup_{k}\eps_{n_{k}}(r)\leq\limsup_{k}\eps_{n_{k}}([r])
     =\eps_{\cu}([r])\to0,\quad r\to\infty.$$
 \end{proof}
 
 We observe here that, by making use of the function $\eps$ considered
 above, it is possible to construct complete metrics on the family
 of equivalence classes of Lip-normed unital \cst{s}.
 
 \begin{Dfn}
     Let $\ca$, $\cb$ be Lip-normed unital \cst{s}, with $\eps$-functions 
     $\eps_{\ca}$, $\eps_{\cb}$, and set
     $$
     \dist^{\eps}_{p}(\ca,\cb):=\max\{\dist_{p}(\ca,\cb), \|\eps_{\ca} -
     \eps_{\cb}\|\},
     $$
     where the norm is the sup norm.
 \end{Dfn}
 
 \begin{Cor}\label{Cor:completeness}
     $\dist^{\eps}_{p}$, $p\geq2$, is a complete metric on the
     family of equivalence classes of Lip-normed unital \cst{s}.
 \end{Cor}
 
 \begin{proof}
     The properties of a metric are obviously satisfied.  Given a
     sequence $\ca_{n}$ of Lip-normed unital \cst{s}, Cauchy w.r.t.
     $\dist_{p}^{\eps}$, the corresponding sequence $\eps_{n}$ is
     uniformly convergent, hence condition $(iii)$ of Corollary
     \ref{Cor:inherited} is satisfied, implying that
     $\ell^{\infty}_{L}(\ca_{n},\cu)$ is a $C^{*}$-algebra.  By
     Proposition \ref{C*ultra} we get the thesis.
 \end{proof}
 
\subsection{Counterexamples}\label{examples}

This section is mainly devoted to the proof of Theorem \ref{noncompl} 
via suitable counterexamples. Also, examples showing the 
non-equivalence of the $f$-Leibniz condition with the $\eps(r)\to0$ 
condition are given.

\subsubsection{Example 1}\label{ex1}

We give here an example of a Cauchy sequence of Lip-normed unital
C$^{*}$-algebras w.r.t. the complete quantum Gromov-Hausdorff distance
$\dist_{\infty}$ which does not converge to a C$^{*}$-algebra.

Let us denote by $\cc$ the algebra of $2\times2$ matrices, and by
$\cc_{0}$ the subspace of $\cc$ consisting of all matrices whose diagonal
part is a multiple of the identity.  Then we let $\cb$ be a
C$^{*}$-algebra acting faithfully on a Hilbert space $\ck$, and denote
by $\ca$ the C$^{*}$-algebra $\cc\otimes\cb$, acting on
$\ch:=\bc^{2}\otimes\ck$, and by $\ca_{0}$ the subspace of $\ca$ given
by $\cc_{0}\otimes\cb$.

Let us now assume that $\cb$ is Lip-normed, with Lip-seminorm $L$, and
define on $\ca$ the functionals
\begin{align*}
    \left\|\begin{pmatrix}a&b \\ c&d\end{pmatrix}\right\|_{n} & :=
\max\left\{ \left\|\frac{a+d}{2}\right\|_{L},n\left\|\frac{a-d}{2}\right\|_{L},
\|b\|_{L},\|c\|_{L}\right\},\quad a,b,c,d\in\cb
\\
L^{n}\left(\begin{pmatrix}a&b \\ c&d\end{pmatrix} \right) & :=
\inf_{\l\in\br} \left\|\begin{pmatrix}a-\l &b \\
c&d-\l\end{pmatrix}\right\|_{n},\quad a,b,c,d\in\cb
\end{align*}

Let us remark that in the following, besides the trivial case $\cb=\bc
I$, we shall consider the case in which $\cb$ is UHF (cf.  Remark
\ref{rem:UHF}).  The existence of a Lip-seminorm on such algebras has been
proved in \cite{AC1}.

\begin{Lemma} 
    $L^{n}$ is a Lip-seminorm on $\ca$.  All these seminorms coincide
    on $\ca_{0}$.
\end{Lemma}

\begin{proof} 
    Obvious.
\end{proof}
    
\begin{Thm}\label{mainbis}
    The sequence $(\ca,L^{n})$ converges in the complete quantum 
    Gromov-Hausdorff distance $\dist_{\infty}$ to $(\ca_{0},L^{1})$.
\end{Thm}

\begin{proof}
    Let us consider the seminorms $\tilde{L}^{n}$ on $\ca_{0}\oplus\ca$:
    \[
    \tilde{L}^{n}(A_{0}\oplus A)=
    \max\{ L^{1}(A_{0}), L^{n}(A),n\|A-A_{0}\|_{1} \},\quad 
    A_{0}\in\ca_{0},A\in\ca.
    \]
    Clearly 
    \[
    \min_{A_{0}\in\ca_{0}}\tilde{L}^{n}(A_{0}\oplus A)= L^{n}(A),\quad
    \min_{A\in\ca}\tilde{L}^{n}(A_{0}\oplus A)= L^{1}(A_{0}),
    \]
    where the first minimum is attained for $A_{0}=
    \begin{pmatrix}(a_{11}+a_{22})/2&a_{12} \\ a_{21}&(a_{11}+a_{22})/2
    \end{pmatrix}$, and the second minimum is attained for $A=A_{0}$.
    This means that $\tilde{L}^{n}$ induces $L^{1}$ on $\ca_{0}$ 
    and $L^{n}$ on $\ca$.
    
    Since $\ca_{0}\subset\ca$, $UCP_{p}(\ca)$ projects onto 
    $UCP_{p}(\ca_{0})$, the projection being simply the restriction 
    to $\ca_{0}$: $\f_{0}:=\f|_{\ca_{0}}$, $\f\in UCP_{p}(\ca)$. 
    Therefore, the distance between $UCP_{p}(\ca_{0})$ and $UCP_{p}(\ca)$ 
    induced by $\tilde{L}^{n}$ is majorised by the supremum, on 
    $\f\in UCP_{p}(\ca)$, of the distance between $\f$ and 
    $\f_{0}=\f|_{\ca_{0}}$. 
    Now
    \begin{align*}
	\r_{\tilde{L}^{n}}(\f_{0}\oplus 0,0\oplus\f)
	&=\sup_{\|A_{0}\oplus A\|_{\tilde{L}^{n}}\leq1}
	\|\langle\f,A_{0}-A\rangle\|\\
	&\leq\sup_{\|A_{0}\oplus A\|_{\tilde{L}^{n}}\leq1}\|A_{0}-A\|\\
	&\leq\sup_{\|A_{0}\oplus A\|_{\tilde{L}^{n}}\leq1}c\|A_{0}-A\|_{1}
	\leq\frac{c}{n},
    \end{align*}
    where we may take $c$ equal to the diameter of 
    $S(\cb)$ w.r.t. $L$. This implies that
    $$
    \dist_{\infty}((\ca,L^{n}),(\ca_{0},L^{1}))
    \leq\sup_{p\in\bn}\r_{\tilde{L}^{n}}^{H}(UCP_{p}(\ca_{0}),UCP_{p}(\ca))
    \leq \frac{c}{n},
    $$
    i.e. the thesis.
\end{proof}

We prove now that $\ca_{0}$ is not a C$^{*}$-algebra up to complete
order isomorphism.  To do so, we need the notion of injective envelope
for operator systems, due to Hamana \cite{Ha}

\begin{Thm}\label{main}
    $\ca_{0}$ is not completely order isomorphic to a C$^{*}$-algebra.
\end{Thm}

\begin{Lemma}
    The injective envelope of $\ca_{0}$ contains $\ca$.
\end{Lemma}

\begin{proof}
    Let $\pi:\cb(\ch)\to\cai(\ca_{0})$ be a completely positive
    projection, existing by injectivity of $\cai(\ca_{0})$.  We will
    show that $\pi$ is the identity on $\ca$.  Choose $b\in\cb_{+}$
    and a unit vector $\xi$ in the Hilbert space $\ck$.  If $u$
    denotes the injection of $\bc\to\ck$ such that $\l\mapsto\l\xi$,
    we may construct the map $\f:\cb(\ch)\to \cc$ given by
    $$
    \f(a)=\begin{pmatrix}u^{*}&0 \\ 0&u^{*}\end{pmatrix}a
    \begin{pmatrix}u&0 \\ 0&u\end{pmatrix}.
    $$
    Let us observe that $\f$ is completely positive and that when $a$
    is written as a $\cb$-valued $2\times2$ matrix we have
    $$
    \f\begin{pmatrix}a_{11}&a_{12} \\ a_{21}&a_{22}\end{pmatrix}
    =\begin{pmatrix}(\xi,a_{11}\xi)&(\xi,a_{12}\xi) \\ 
    (\xi,a_{21}\xi)&(\xi,a_{22}\xi)\end{pmatrix}.
    $$
    We then consider the map $\psi:A\in\cc\to \psi(A)\in\cc$ given by 
    $\psi(A)=\f(\pi(A\otimes b))$, and notice that $\psi$ is 
    completely positive and, when $A\in \cc_{0}$, we have $\pi(A\otimes 
    b)=A\otimes b$, hence 
    \begin{equation}\label{psirelation}
	\psi(A)=(\xi,b\xi)A.
    \end{equation}
    Let us show that this relation holds for any $A\in \cc$.  Indeed
    this is clearly true when $(\xi,b\xi)=0$, since a positive map
    vanishing on the identity is zero.  When $(\xi,b\xi)\ne 0$, the
    map $\frac1{(\xi,b\xi)}\psi$ is a completely positive map from
    $\cc$ to $\cc$ which is the identity on $\cc_{0}$ and, since the
    injective envelope of $\cc_{0}$ is $\cc$, it has to be the
    identity anywhere.  A simple calculation shows that relation
    (\ref{psirelation}) may be rewritten as $(\xi,(\pi(A\otimes
    b)_{ij}-a_{ij}b) \xi)=0$, $i,j=1,2$.  By the arbitrariness of
    $\xi$ we get $\pi(A\otimes b)=A\otimes b$, and by the
    arbitrariness of $b\in\cb$ we get the thesis.
\end{proof}

\begin{proof}[Proof of Theorem \ref{main}]
    Let us recall Proposition 15.10 in \cite{Paulsen}: given an
    inclusion $\cb\subseteq S\subset\cb(\ch)$, where $\cb$ is a unital
    C$^{*}$-algebra and $S$ is an operator system, then $\cb$ is a
    subalgebra of $\cai(S)$.  This implies that if $S$ is an operator
    system that can be represented as a unital C$^{*}$-algebra $\cb$
    acting on $\ch$, the immersion of $\cb$ in $\cai(\cb)$ is a
    $^{*}$-monomorphism, namely the product structure of $S$ making it
    a C$^{*}$-algebra is the one given by the immersion in its
    injective envelope.
    
    Then, posing $S=\cai(\ca_{0})$ and $\cb=\ca$ in the same Proposition,
    one gets that the product on $\ca_{0}$ given by the immersion in
    $\cai(\ca_{0})$ is the same as that given by the immersion in
    $\ca$, namely $\ca_{0}$ is not a subalgebra of its injective
    envelope.  By the argument above, it is not an algebra.
\end{proof}

\begin{Cor}
    The space of equivalence classes of C$^{*}$-algebras endowed with the 
    metric $\dist_{\infty}$ is not complete.
\end{Cor}

\begin{rem}\label{rem:UHF}
    The preceding example works well also in the case $\cb=\bc$. 
    However, in the finite-dimensional case, the replacement of the
    distance between state spaces with the distance between (the
    closure of) pure states, like the distance $\dist^{e}_{q}$
    considered by Rieffel in \cite{Rieffel00} after Proposition 4.9,
    would destroy the example, since the sequence is not Cauchy w.r.t.
    such distance.  One could therefore think that, endowing
    C$^{*}$-algebras with the appropriate distance, completeness may
    follow.  But this is not true, since, choosing $\cb$ as a UHF
    algebra, the pure states are dense, namely the mentioned
    replacement would have no effect.
    
    Let us also mention that when $\cb=\bc$, namely $\ca_{0}=\cc_{0}$,
    such operator system is not even order isomorphic to a
    C$^{*}$-algebra.  Indeed its state space is two dimensional and
    has the convex structure of a disc, while the only C$^{*}$-algebra
    with two-dimensional state-space is $\bc^{3}$, whose state space
    has the convex structure of a triangle.  This means that even
    replacing $\dist_{\infty}$ with $\dist_{p}$ the set of Lip-normed
    unital C$^{*}$-algebras is still non-complete.
\end{rem}

\subsubsection{Example 2}\label{ex2}
We give here an example of a Cauchy sequence of Lip-normed unital
C$^{*}$-algebras w.r.t. the complete quantum Gromov-Hausdorff distance
$\dist_{\infty}$ which converges to a C$^{*}$-algebra, but the
C$^{*}$-structure is not inherited.

The sequence $\{\ca_{n}\}_{n\in\bn}$ is made of the constant algebra
$\bc^{3}$ endowed with the following seminorms:
$$
L_{n}(a,b,c)=\|\frac{a-b}{2},n(\frac{a+b}{2}-c)\|_{2},
$$
where $\|\cdot\|_{2}$ is the Euclidean norm.  It is not difficult to
show that the sequence converges, in any $\dist_{p}$, to the Lip-normed
unital C$^{*}$-algebra $\ca_{\infty}$ consisting of $\bc^{2}$ with the
seminorm $L_{\infty}(\a,\b) = |\frac{\a-\b}{2}|$.  Indeed, let us
consider on $\bc^{3}\oplus \bc^{2}$ the seminorm
$$
\tilde{L}_{n}(a,b,c,\a,\b)=
\max \{ L_{n}(a,b,c), L_{\infty}(\a,\b),
n|a-\a|,n|b-\b|,n|c-\frac{\a+\b}{2}| \}.
$$
Clearly $\tilde{L}_{n}$ induces $L_{n}$ on $\ca_{n}$ and  
$L_{\infty}$ on $\ca_{\infty}$ and, reasoning as in the previous example, we 
get $\dist_{\infty}(\ca_{n},\ca_{\infty})\leq\frac1n$.

Now we compute the ultraproducts.  Since we have a sequence of
finite-dimensional constant spaces, for any free ultrafilter $\cu$, the
ultraproduct coincides with $\bc^{3}$, where we can represent any
element with the constant sequence \cite{AkKh}.  Then the Lip-ultraproduct
consists of those sequences constantly equal to $(a,b,c)$ for which
$L_{n}(a,b,c)$ is bounded, i.e. $c=\frac{a+b}{2}$.  Therefore, setting
$$
\ca_{0}=\{(\a,\b,\frac{\a+\b}{2})\in\bc^{3} : \a,\b\in\bc\},
$$ 
the inclusion of the Lip-ultraproduct in the ultraproduct is
given by $\ca_{0}\subset\bc^{3}$, for any free ultrafilter $\cu$.  Since
$\ca_{0}$ is not a subalgebra of $\bc^{3}$, the limit does not
inherit the C$^{*}$-structure.

Let us remark that the previous results are not in contradiction,
since the map $(a,b)\in\bc^{2} \mapsto (a,b,(a+b)/2)$ is clearly a
complete order isomorphism, namely $\ca_{0}$ and $\ca_{\infty}$ are
completely order isomorphic.

\begin{rem}
    The previous example consists of abelian C$^{*}$-algebras
    converging to an abelian C$^{*}$-algebra, therefore one could
    expect it corresponds to the Gromov-Hausdorff convergence of the
    spectra.  But if this were true the Lip-ultraproduct would
    correspond to the ultralimit, hence would be a C$^{*}$-algebra in
    a natural way.  This apparent contradiction is due to the fact
    that the approximating state spaces (triangles) converge to the
    limit state space (segment) like a triangle flattening on its
    base, namely the upper vertex converges to the middle point of the
    basis.  Therefore the spectra do not converge Gromov-Hausdorff.
\end{rem}

\subsubsection{Example 3}\label{ex3}
 
We conclude with an example of a converging sequence of
C$^{*}$-algebras where the limit inherits the C$^{*}$-structure,
however no $f$-Leibniz condition is satisfied, namely there is no
function $f$ such that all algebras satisfy the same $f$-Leibniz
condition.  According to Proposition \ref{f-Leibniz}, it is sufficient
to exhibit a converging sequence for which the functions $\eps_{n}$
are eventually zero, but converge pointwise to a nowhere zero function
infinitesimal at $+\infty$.
 
As in the previous examples, the sequence will consist of a constant 
algebra with varying Lip-seminorms.

The C$^{*}$-algebra $\ca$ is made of sequences $A=\{A_{k}\}_{k\in\bn}$ of
$2\times 2$ matrices converging to a multiple of the identity.

On the C$^{*}$-algebra $\ca$ let us consider the (possibly infinite) functionals 
\begin{align*}
    \|A\| & = \sup_{k} \|A_{k}\|, \\
    \tn{A}&=\sup_{k}k \|A_{k}\|,\\
    L(A)&=\min_{\l\in\bc} \tn{A-\l I}.
\end{align*}
and the dense subspace $\ca_{0}$ of the elements for which $L(A)<\infty$.  

Let us observe that if $\tn{A-\a I}<\infty$ then $A_{k}\to \a I$,
hence 
\begin{equation}\label{A-alpha}
    |\a|=\lim_{k}\|A_{k}\|\leq \sup_{k}\|A_{k}\|=\|A\|.
\end{equation}

\begin{Lemma}
   $L$ is a Lip-seminorm, and satisfies the inequality
    \begin{equation}\label{Leibniz}
	L(AB)\leq L(A)\|B\|+\|A\|L(B).
    \end{equation}
\end{Lemma}

\begin{proof}
    Clearly $\tn{\cdot}$ is a lower semicontinuous norm on $\ca_{0}$,
    hence $L$ is a lower semicontinuous seminorm vanishing only on the
    multiples of the identity.  Let us observe that $\cb:=
    \{B:\tn{B}\leq 1\}$ is the image of the unit ball under the
    compact operator sending $\{A_{k}\} \mapsto \{\frac1kA_{k}\}$,
    hence it is totally bounded.  Consider $\{A:L(A)\leq1,\|A\|\leq
    1\}$.  Then $\tn{A-\a I}\leq1$ for a suitable $\a$.  Making use of
    inequality (\ref{A-alpha}), we get $A\in \cup_{|\a|\leq 1} (\a
    I+\cb)$, showing that such set is totally bounded, i.e. $L$ is a
    Lip-seminorm.
    
    Concerning inequality (\ref{Leibniz}), we have, for $A,B\in 
    \ca_{0}$ with $\tn{A-\a}=L(A)$, $\tn{B-\b}=L(B)$, 
    \begin{align*}
	L(AB)&\leq\tn{AB-\a\b}=\tn{(A-\a)B+\a(B-\b)}\\
	&\leq L(A)\|B\|+|\a|L(B)\leq L(A)\|B\|+\|A\|L(B),
    \end{align*}
    where we used inequality (\ref{A-alpha}).
\end{proof}

Now we consider a new sequence of Lip-seminorms on $\ca$:
$$
L_{n}(A)= \max \{ L(A), \sup_{k< n}\ell_{k}(A_{k}) \},\quad n\in\bn\cup 
\{\infty\},
$$
where 
$$
\ell_{k}\begin{pmatrix}a&b \\ c&d\end{pmatrix}=k^{3}|a-d|.
$$
Clearly each $L_{n}$ is again a Lip-seminorm, and, for finite $n$, it
still satisfies an $f$-Leibniz condition (cf.  Proposition
\ref{f-Leibniz}), being a finite rank perturbation of $L$.  

In the following we shall denote by $\ca_{n}$ the Lip-normed unital
C$^{*}$-algebra $(\ca, L_{n})$, $n\in\bn\cup\{\infty\}$.

First we observe that, for any free ultrafilter $\cu$, we may identify
the Lip-ultraproduct $\ell_{L}^{\infty}(\ca_{n},\cu)$ with
$\ca_{\infty}$.  Indeed, given $\{A^{n}\}_{n\in\bn}\subset \ca$ with
$\|A^{n}\|\leq 1$ and $L_{n}(A^{n})\leq1$, we have shown that it lies
in a compact set, namely $\lim_{\cu}A^{n}$ exists, and we call it $A$. 
We can therefore identify the class of the sequence $\{A^{n}\}$ in
$\ell_{L}^{\infty}(\ca_{n},\cu)$ with the class of the sequence
constantly equal to $A$.  As a consequence the C$^{*}$-structure is
inherited.

Now we show that indeed $\{\ca_{n}\}$ converges in the complete
quantum Gromov-Hausdorff distance $\dist_{\infty}$ to $\ca_{\infty}$. 
Take on $\ca\oplus\ca$ the seminorm
$$
\tilde{L}_{n}(A,B) = \max\{ L_{n}(A),L_{\infty}(B), n\|A-B\| \},
$$
which is clearly a Lip-seminorm.  It is easy to see that it induces
$L_{n}$ on the first summand, the minimum being attained for
$B_{k}=A_{k}$, $k\leq n$, $B_{k}=\a I$, $k> n$, with $L(A)=\tn{A-\a
I}$.  Analogously, it induces $L_{\infty}$ on the second summand.

As in the first example, we get 
$$
\r_{\tilde{L}^{n}}(\f\oplus 0,0\oplus\f)
\leq\sup_{\|A\oplus B\|_{\tilde{L}^{n}}\leq1}\|A-B\|\leq\frac{1}{n},
$$
hence
$$
\dist_{\infty}(\ca_{n},\ca_{\infty})
\leq\sup_{p\in\bn}\r_{\tilde{L}^{n}}^{H}(UCP_{p}(\ca\oplus0),UCP_{p}(0\oplus
\ca)) \leq \frac{1}{n},
$$
i.e. the thesis.

It only remains to show that $(\ca, L_{\infty})$ does not satisfy any
$f$-Leibniz condition, i.e. by Proposition \ref{f-Leibniz}, that we
can find an element $A$ with finite Lip-seminorm such that
$L_{\infty}(A^{*}A)$ is infinite.  Taking $A=\{A_{k}\}$,
$A_{k}=\begin{pmatrix}0&1/k \\ 0&0\end{pmatrix}$, we have
$L_{\infty}(A)=L(A)=1$, but $L_{\infty}(A^{*}A)=\infty$.

\begin{ack} We would like to thank Marc Rieffel for discussions and
    David Kerr and Hanfeng Li for comments and suggestions.
\end{ack}

%%%%%%REFERENCES%%%%%%%%

\end{document}